\documentclass[12pt, leqno]{amsart}

\usepackage[OT2,T1]{fontenc}
\usepackage{indentfirst}
\usepackage{amstext}
\usepackage{amsthm}
\usepackage{amsopn}
\usepackage{amsfonts}
\usepackage{amsmath}
\usepackage{latexsym}
\usepackage[francais,english]{babel}
\usepackage{amscd}
\usepackage{amssymb}
\usepackage{amsmath}
\usepackage[all,cmtip]{xy}
\usepackage{leftidx}
\usepackage{graphicx}
\usepackage{etoolbox}
\patchcmd{\section}{\normalfont\scshape\centering}{\normalfont\bfseries}{}{}
\patchcmd{\subsection}{-.5em}{.5em}{}{}

\newtheorem{theo}{{Theorem}}[section]
\newtheorem{coro}[theo]{{Corollary}}
\newtheorem{lemma}[theo]{{Lemma}}
\newtheorem{prop}[theo]{Proposition}

\theoremstyle{definition}
\newtheorem{remark}[theo]{\textbf{Remark}}
\newtheorem{defn}[theo]{Definition}
\newtheorem{example}[theo]{Example}
\numberwithin{equation}{section}

\newtheorem{notation}[theo]{Notation}
\newtheorem{question}[theo]{Question}

\DeclareFontEncoding{OT2}{}{} 

\begin{document}
\tolerance 400 \pretolerance 200 \selectlanguage{english}

\title{$K3$ surfaces, cyclotomic polynomials and orthogonal groups}
\author{Eva Bayer-Fluckiger}
\date{\today}
\maketitle

\centerline {In memory of Nikolai Vavilov}

\begin{abstract} Let $X$ be a complex projective $K3$ surface, and let $T_X$ be its transcendental lattice; the characteristic polynomials of the isometries of $T_X$ induced
by automorphisms of $X$ are powers of  cyclotomic polynomials. Which powers of cyclotomic polynomials occur ? The aim of this note is to answer this question, as well
as related ones, and give an alternative approach to some results of Kondo, Machida, Oguiso, Vorontsov, Xiao and Zhang; this leads to questions
and results concerning orthogonal groups of lattices.

\medskip

\end{abstract}

\small{} \normalsize

\medskip

\selectlanguage{english}
\section{Introduction}

If $X$ is a projective $K3$ surface over the complex numbers; we denote by $S_X$ its Picard lattice and by $T_X$ its transcendental lattice; if $a : X \to X$
is an automorphism, then $a$ induces an isometry $a^*$ of the lattice $H^2(X,{\bf Z})$, and the characteristic polynomial of the restriction of $a^*$ to
$T_X$ is a power of a cyclotomic polynomial (see Proposition \ref{cyclotomic minimal polynomial}).

\medskip
Let $m, r$ be integers with $m \geqslant 3$ and $r \geqslant 1$, and let $C = \Phi_m^r$ (where $\Phi_m$ is the $m$-th cyclotomic polynomial). 

\medskip
\noindent
{\bf Proposition 1.} {\it Assume that ${\rm deg}(C) \leqslant 20$. Then there exists an automorphism $a : X \to X$ of a projective $K3$ surface $X$ such that the characteristic polynomial of
the restriction of $a^*$ to $T_X$ is equal to $C$.}

\medskip
We denote by ${\rm Aut}(X)$ the group of automorphisms of the $K3$ surface $X$, and by ${\rm Aut}_s(X)$ the subgroup of ${\rm Aut}(X)$ acting trivially on $T_X$. We
have the exact sequence 
$$1 \to {\rm Aut}_s(X) \to {\rm Aut}(X) \to M_X \to 1,$$ where $M_X$ is a finite cyclic group (see Nikulin, \cite{N}, Theorem 10.1.2); we denote by $m_X$ the order of $M_X$. 

\medskip
\noindent
{\bf Corollary 1.} {\it Let $m \geq 4$ be an even integer such that $\varphi(m) \leqslant 20$. Then  there exists a projective $K3$ surface $X$ 
with $m_X = m$.}

\medskip It is well-known that there exist $K3$ surfaces $X$ with $m_X = 1,2$ (see for instance \cite{H}, Corollary 15.2.12), but as far as I know, the following question is open :

\medskip
\noindent
{\bf Question 1.} {\it Let $m > 1$ be an odd integer. Does there exist a projective $K3$ surface $X$ 
with $m_X = m$ ?}

\medskip
In \cite{MO}, Machida and Oguiso obtain several results on related topics; see Remark \ref{finite order MO}, Proposition \ref{no 60} and Remark \ref {OM list} for details. 

\bigskip Following Vorontsov \cite{V} and Kondo \cite{Ko}, we consider automorphisms that act trivially on the Picard lattice. 
Let $N_X$ be  the kernel of ${\rm Aut}(X) \to {\rm O}(S_X)$; this is a finite cyclic group, that can be identified with a subgroup of $M_X$; we
denote by $n_X$ the order of $N_X$. 
The following proposition is proved in 
 \S \ref{Vorontsov} :

\medskip
\noindent
{\bf Proposition 2.} {\it There exists an automorphism $a : X \to X$ of a  projective $K3$ surface $X$ such that $a^*$ is the identity on $S_X$ and that the characteristic polynomial of
the restriction of $a^*$ to $T_X$ is equal to $C$ if and only if the following conditions hold

\medskip
{\rm (i)} $C(-1)$ is a square.

\medskip
{\rm (ii)} If $C(1) = 1$, then  ${\rm deg}(C)  \equiv \ 4  \ {\rm (mod \ 8)}$. }

\medskip
The possible values of $n_X$ can be deduced from Proposition 2,  extending 
previous results of Vorontsov \cite{V}, Kondo \cite{Ko} and Oguiso-Zhang \cite{OZ}; see \S \ref{cyclic}. Note that $n_X$ divides $m_X$, since $N_X$ can be identified
with a subgroup of $M_X$. This suggests the following question

\medskip
\noindent
{\bf Question 2.} What are the possible values of the pairs $(m_X,n_X)$ ?

 \bigskip
 The proofs of the above propositions use some arithmetic results (see below), as well as
 the surjectivity of the period map, the strong Torelli theorem, and some results of McMullen \cite{Mc3}.
 
 \medskip
 The arithmetic results are valid in a greater generality than the one needed for the applications to $K3$ surfaces. For instance,
 in order to prove Proposition 2, we introduce the following property, called property (P 2)  :
 
 \medskip
 Let $R,S \geqslant 0$ be integers  such that  such that $R  \equiv S \ {\rm (mod \ 8)}$, and set $N = R+S$; suppose that  
 ${\rm deg}(C) < N$.

 \medskip 
 Let $c \geqslant 0$ be an even integer with $c \leqslant {\rm deg}(C)$ such that $c \leqslant R$ and ${\rm deg}(C) - c \leqslant S$.

 \medskip 
 
 \noindent
 {\bf Definition.} We say that {\it property {\rm (P 2)} holds} if there exists an even, unimodular lattice $L$ 
 of signature $(R,S)$ and an isometry $t : L \to L$ such that
 
 \medskip
 $\bullet$ The characteristic polynomial of $t$ is $C(X)(X-1)^{N - {\rm deg}(C)}$.
 
 \medskip
 $\bullet$ The signature of the sublattice ${\rm Ker}(C(t))$ is $(c,{\rm deg}(C)-c)$. 
 
 \bigskip In the application to $K3$ surfaces, we have $R = 3$, $S = 19$ and $c = 2$, and
 Proposition 2 is a consequence of the following result : 
 
 \bigskip
 \noindent
 {\bf Theorem 2.} {\it Property {\rm (P 2)}  holds if and only if the following conditions are satisfied : 

\medskip
{\rm (i)} $C(-1)$ is a square.

\medskip
{\rm (ii)} If $C(1) = 1$, then  ${\rm deg}(C)    \equiv \ 2c  \ {\rm (mod \ 8)}$. }

\bigskip
In particular, Property 2 always holds if $C(-1)$ is a square and  $C(1) > 1$. 

\medskip
If $m$ is an odd prime number, this can be deduced from a result of Brandhorst and Cattaneo, \cite{BC}, Theorem 1.1.
Note that Theorem 2 gives a partial answer to a question of this paper (see  \cite{BC}, Outlook).

\medskip
With the same notation, we introduce property (P 1) :

\medskip
 
 \noindent
 {\bf Definition.} We say that {\it property {\rm (P 1)} holds} if there exists an even, unimodular lattice $L$ 
 of signature $(R,S)$ and an isometry $t : L \to L$ such that
 
 \medskip
 $\bullet$ The characteristic polynomial of $t$ is divisible by $C$. 
 
 \medskip
 $\bullet$ The signature of the sublattice ${\rm Ker}(C(t))$ is $(c,{\rm deg}(C)-c)$. 
 
 \medskip
 \noindent
 {\bf Theorem 1.} If $C(1)C(-1) > 1$, then Property (P 1) holds. 
 
 \medskip If $C(1) = C(-1) = 1$ and
 $R = 0$ or $S = 0$, then Property (P 1) does not always hold, but
 the {\it indefinite} case seems to be open :
 
 \medskip
 \noindent
 {\bf Question 3.} Does Property (P 1) always hold when $R > 0$ and $S > 0$ ?
 
 \medskip
  A modified version
 of this property is used to prove Proposition 1. A more tractable question is to ask for isometries of {\it finite order}; this
 leads to the following definition :
 
 \medskip
 
  \noindent
 {\bf Definition.} We say that {\it property {\rm (P 1')} holds} if there exists an even, unimodular lattice $L$ 
 of signature $(R,S)$ and an isometry $t : L \to L$ of finite order such that
 
 \medskip
 $\bullet$ The characteristic polynomial of $t$ is divisible by $C$. 
 
 \medskip
 $\bullet$ The signature of the sublattice ${\rm Ker}(C(t))$ is $(c,{\rm deg}(C)-c)$. 
 
 \medskip Note that Properties (P 1) and (P 1') are equivalent if 
 $R = 0$ or $S = 0$,  since then all isometries are of finite order. If
 $R > 0$ and $S > 0$,  it is possible that Property (P 1) always holds - however, this is not the case for
 Property (P 1'), as shown by the following example :
 
 \medskip
 \noindent
 {\bf Example.} Let $C = \Phi_{60}$, and set $R = 3$, $S = 19$ and $c = 2$. Then
 
 \medskip
 $\bullet$ Property (P 1) holds (see Example \ref{for 60 cyclotomic}),
 
 \medskip
 $\bullet$ Property (P 1') does not hold (see Proposition \ref{60.1}). 
 
 \medskip
 Note that the fact that property (P 1') does not hold implies the well-known fact that there does not exist any projective $K3$ surface $X$
 having an automorphism $a \in {\rm Aut}(X)$ of {\it finite order} such that $a^*$ induces multiplication by a primitive $60$-th root of unity on $T_X$; 
 see Machida and Oguison \cite {MO}, Xiao \cite{X}, and  Zhang \cite{Z}; see also Proposition \ref{no 60}.

 \medskip

 It can be useful to replace ``of finite order'' by  ``of {\it order $m$}''. This point of view
 is taken by several authors; see the paper of Brandhorst \cite{Bra} and the references therein; see also \S \ref{P 1 section}, \S \ref{0 section}
 and \S \ref{2 section}. 
 
 

 \medskip
 The paper is organized as follows. The first two sections are mainly preliminary, recalling notions and results on $K3$ surfaces and
 isometries of lattices; \S \ref{isometries} and \S \ref{Salem section} also recall some results of \cite{B 21} and \cite{B 22}, and give
 some examples that are used in the paper. Theorem 2 is proved in \S \ref{P 2} (see Theorem \ref {Theorem 2}), and Proposition 2
 in \S \ref{Vorontsov} (see Proposition \ref{Kondo}). A stronger form of Theorem 1 is in \S \ref{P 1 section}, see Theorem \ref{theorem 1}.
 Proposition 1 is proved in \S \ref{Prop 1} and \S \ref{Prop 1 2p}. The last section concerns the possible values of $m_X$ and $n_X$,
 a discussion of some results of Kondo and Vorontsov, as well as a generalization of these results, and some open questions.
 
 \medskip The proofs use  results of \cite{B 21}, \cite{B 22}, of McMullen, \cite{Mc3} and  of Takada \cite{T}.
 
 
 

\section{$K3$ surfaces}\label{K3}

We recall some notation and basic facts on $K3$ surfaces and their automorphisms; see
\cite {H} and  \cite {K}
for details. 

\medskip A $K3$ surface $X$ is a simply-connected compact complex surface with trivial
canonical bundle. We
have the Hodge decomposition $$H^2(X,{\bf C}) = H^{2,0}(X) \oplus H^{1,1}(X) \oplus H^{0,2}(X)$$ with ${\rm dim} \ H^{2,0} = {\rm dim} \ H^{0,2} = 1$
and ${\rm dim} \ H^{1,1} = 20$.
The {\it Picard lattice} of 
$X$ is by definition $$S_X = H^2(X,{\bf Z}) \cap H^{1,1}(X).$$
The intersection form $H^2(X,{\bf Z}) \times H^2(X,{\bf Z}) \to {\bf Z}$ of $X$
is an even unimodular lattice of signature $(3,19)$. Such a form is unique up to isomorphism (see for instance \cite{S}, Chap. V, th. 5).
The {\it transcendental lattice} $T_X$ is by definition the primitive sublattice of $H^2(X,{\bf Z})$ of minimal rank such that $T_X \otimes_{\bf Z} {\bf C}$
contains $H^{2,0}(X) \oplus H^{0,2}(X)$. 
Assume that $X$ is {\it projective}; then lattices $S_X$ and $T_X$ are orthogonal to each other, and the orthogonal sum $S_X \oplus T_X$ is of finite
index in $H^2(X,{\bf C})$, the signature of  $S_X$ is $(1,\rho_X-1)$, and the signature of $T_X$ is $(2,20-\rho_X)$, where $\rho_X$ is
the rank of $S_X$. 

\medskip
If $a : X \to X$ is an automorphism, then $a^* : H^2(X,{\bf C}) \to H^2(X,{\bf C})$ respects the Hodge decomposition and is an isometry of the intersection form; hence $a^*$ is also an isometry of the
lattices $S_X$ and $T_X$. 

\medskip The following is a result of Oguiso, \cite{O 8}, Theorem 2.4.

\begin{prop}\label{cyclotomic minimal polynomial} Let $a : X \to X$ be an automorphism of a projective $K3$ surface. Then the
minimal polynomial of the restriction of $a^*$ to $T_X$ is a cyclotomic polynomial. 

\end{prop}

\noindent
{\bf Proof.} The minimal polynomial of the restriction of $a^*$ to $T_X$ is irreducible (see Oguiso \cite{O 8}, Theorem 2.4, (1)).
Since $X$ is projective, $a^*|T_X$ is of finite order (cf. Nikulin \cite{N}, Theorem 10.1.2,  \cite{O 8}, Theorem 2.4 (4), or \cite{H}, Corollary 3.3.4 or Corollary 15.1.10), hence 
its minimal polynomial is a cyclotomic polynomial.

\section{Isometries of lattices}\label{isometries}

In this section we summarize some notions and results from \cite{GM}, \cite{B 21} and \cite{B 22} in the special cases needed in this paper.

\medskip

A {\it lattice} is a pair $(L,q)$, where $L$ is a free ${\bf Z}$-module of finite rank, and $q : L \times L \to {\bf Z}$
is a symmetric bilinear form; it is {\it unimodular} if ${\rm det}(q) = \pm 1$, and {\it even} if $q(x,x)$ is an
even integer for all $x \in L$. The following lemma is well-known :

\begin{lemma}\label{mod 8} Let $L$ be an even unimodular lattice of signature $(r,s)$. Then $r  \equiv s \ {\rm (mod \ 8)}$.

\end{lemma}

\noindent
{\bf Proof.} See for instance \cite{S}, Chap. V, cor. 1.

\medskip
Let $n \geqslant 1$ be an integer and let $F \in {\bf Z}[x]$ be a monic polynomial of degree $2n$ such that $F(x) = x^{2n}F(x^{-1})$. We say that
$F$ {\it satisfies condition} (C 1) if $|F(1)|$, $|F(-1)|$ and $(-1)^nF(1)F(-1)$ are squares. 

\begin{lemma}\label{condition c1} {\rm (Gross-McMullen)} Let $L$ be an even unimodular lattice and let $t : L \to L$ be an isometry  with characteristic
polynomial $F$. Then $F$ satisfies condition {\rm (C 1)}. 

\end{lemma} 

\noindent
{\bf Proof.} See \cite{GM}, Theorem 6.1 or \cite {B 22}, Corollary 2.3. 

\medskip Assume that $F$ is a {\it product of cyclotomic polynomials}, and that 
$F(1) \not = 0$ or $F(-1) \not = 0$.  Let us write $F = F_1 F_0$, where $F_1 \in {\bf Z}[x]$ is a monic polynomial such that $F_1(1)F_1(-1) \not = 0$, and $F_0(x) = (x-1)^{n_+}$ or $F_0(x) = (x+1)^{n_-}$, where $n_+,n_- \geqslant 0$ are integers. Set $D_0 = (-1)^nF_1(1)F_1(-1)$.  Let $I$ be the set of irreducible factors of $F$ over $\bf Q$. 

\medskip
Following \cite{B 22}, we associate to $F$ a finite group $G_F$; we start by defining a set $\Pi_{f,g}$ for all $f,g \in I$, as follows. We say that a monic polynomial $h$ is ($\pm$)-{\it symmetric} if
$h(x) = \pm x^{{\rm deg}(h)}h(x^{-1})$. We also use the terminology {\it symmetric} for (+)-symmetric. 

\medskip
If $f,g \in I$ are such that ${\rm deg}(f) \geqslant 2$, ${\rm deg}(g) \geqslant 2$, then $\Pi_{f,g}$ is the set of prime numbers $p$ such that 
$f  \ {\rm (mod \ {\it p})}$ and $g  \ {\rm (mod \ {\it p})}$
 have a common irreducible ($\pm$)-symmetric 
factor in $ {\bf F}_p[x]$.

\medskip If $f \in I$ is such that ${\rm deg}(f) \geqslant 2$, then $\Pi_{f,x-1}$ is the set of prime numbers $p$ such that $f  \ {\rm (mod \ {\it p})}$
is divisible by $x-1$ 
 in ${\bf F}_p[x]$, and that if $n^+ = 2$, then $D_0 \not = -1$ in ${\bf Q}_p^{\times}/{\bf Q}_p^{\times 2}$. 
 
 \medskip If $f \in I$ is such that ${\rm deg}(f) \geqslant 2$, then $\Pi_{f,x+1}$ is the set of prime numbers $p$ such that $f  \ {\rm (mod \ {\it p})}$
is divisible by $x+1$ 
 in ${\bf F}_p[x]$, and that if $n^- = 2$, then $D_0 \not = -1$ in ${\bf Q}_p^{\times}/{\bf Q}_p^{\times 2}$. 
 
 \medskip

 Let $C(I)$ be the set of maps $c : I \to {\bf Z}/2{\bf Z}$,
 let $C_0(I)$ the set of $c \in C(I)$ such that $c(f) = c(g)$ if $\Pi_{f,g} \not = \varnothing$; note that $C_0(I)$ is an abelian group. Let $G_F$ be the quotient of $C_0(I)$ by the subgroup of constant maps. 
 
 \begin{example} Let $F(x) = \Phi_{15}^2(x) \Phi_3(x) (x-1)^2$. We have $5 \in \Pi_{\Phi_{15},\Phi_3}$ and $3 \in \Pi_{\Phi_3(x),x-1}$, hence $G_F = 0$. 
 
 \end{example} 
 
 When $F$ has no linear factors, then $G_F$ is already defined in \cite{B 21}, and several examples are given in \cite{B 21}, \S 25. Here is another example, that
 will be used in the proof of Proposition \ref{8} :
 
 \begin{example}\label{12 and 60} Let $F = \Phi_{60} \Phi_{12}$. The resultant of $\Phi_{60}$ and $\Phi_{12}$ is $5^4$, and the polynomials $\Phi_{60}$ and
$\Phi_{12}$  ${\rm (mod} \ 5 {\rm )}$  have the common irreducible factors $x^2 + 2x + 4$ and $x^2 + 3x + 4$ in ${\bf F}_5[x]$. These polynomials are not ($\pm$)-symmetric, hence 
$\Pi_{\Phi_{60}, \Phi_{12}} = \varnothing$, and $G_F = {\bf Z}/2{\bf Z}$.

 \end{example}
 
 The following is proved in \cite{B 22}, Corollary 12.4 :
 \begin{theo}\label{basic theorem} Let $r,s \geqslant 0$ be integers such that $r  \equiv s \ {\rm (mod \ 8)}$ and that $r+s = {\rm deg}(F)$. If $G_F = 0$, then there exists an
 even unimodular lattice of signature $(r,s)$ having an isometry with characteristic polynomial $F$. 
 
 \end{theo} 
 
 We need a more precise result : it is not enough to fix the signature of the lattice, we also need information about the {\it signature map} of the isometry. We
 recall this notion from \cite{B 22}, \S 3 and \S 4.

\begin{defn}\label{signature definition} Let $V$ be a finite dimensional vector space over $\bf R$, and let $q : V \times V \to {\bf R}$ be
a non-degenerate quadratic form.
Let $t : V \to V$ be an isometry of $q$. If $f \in {\bf R}[X]$, set $V_f = {\rm Ker}(f(t))$ and let $q_f$ be the restriction of $q$ to $V_f$. 
Let $${\rm sign}_t : {\bf R}[X] \to {\bf N} \times {\bf N}$$
be the map sending $f \in {\bf R}[X]$ to the signature of $(V_f,q_f)$, where $\bf N$ is the set of all nonnegative integers; it is called the {\it signature map
of the isometry} $t$. The signature of $q$ is called the {\it maximum} of the signature map, and the characteristic polynomial of $t$ the
{\it polynomial associated to the signature map}. 

\end{defn}


\begin{example}
 Let $a : X \to X$ be an automorphism of a projective $K3$ surface, and suppose that $a^*|S_X$ is the identity and
that the characteristic polynomial of the restriction of $a^*$ to $T_X$ is $C$. Let $\rho_X$ be the rank of $S_X$, and let $\tau$ be the signature map of $a^*$. Then
we have $\tau(x-1) = (1,\rho_X-1)$ and $\tau(C) = (2,{\rm deg}(C)-2)$. Note that ${\rm deg}(C) = 22 - \rho_X$, hence we also have
$\tau(C) = (2,20-\rho_X)$. 

\end{example}

\begin{theo}\label{signature map theorem}  Let $r,s \geqslant 0$ be integers such that $r  \equiv s \ {\rm (mod \ 8)}$ and that $r+s = {\rm deg}(F)$. Let $\tau$ be a signature
map of maximum $(r,s)$ and associated polynomial $F$.  If $G_F = 0$, then there exists an
 even unimodular lattice  having an isometry with signature map $\tau$.  
 
 \end{theo} 
 
 \noindent
 {\bf Proof.} This is the statement of  \cite{B 22}, Corollary 12.3.

\section{Cyclotomic polynomials and property (P 2)}\label{P 2}

The aim of this section is to prove Theorem 2 of the introduction. We start by recalling some basic properties of
cyclotomic polynomials. Recall that $\Phi_m$ denotes the $m$-th cyclotomic polynomial, and that ${\rm deg}(\Phi_m) = \varphi(m)$. 

\begin{lemma}\label{1 and -1} Let $m$ be an integer with $m \geqslant 3$. We have

\medskip
{\rm (i)} If $m$ is a power of $2$, then $\Phi_m(1) = \Phi_m(-1) = 2$.

\medskip
{\rm (ii)} Let  $p$ be a prime number with $p \not = 2$ and let $k \geqslant 1$ an integer. If $m = p^k$, then $\Phi_m(1) = p$ and $\Phi_m(-1) = 1$; if 
$m = 2p^k$, then $\Phi_m(1) = 1$ and $\Phi_m(-1) = p$.

\medskip
{\rm (iii)} In all other cases, $\Phi_m(1) = \Phi_m(-1) = 1$.

\medskip
{\rm (iv)} If $\Phi_m(1) = \Phi_m(-1) = 1$, then ${\rm deg}(\Phi_m) \equiv 0 \ {\rm (mod \ 4)}$.

\end{lemma} 

\noindent
{\bf Proof.} (i) Let $m = 2^k$ for some integer $k \geqslant 2$. Then $\Phi_m(x) = x^{2^k-1}+1$, hence
$\Phi_m(1) = \Phi_m(-1) = 2$.

\medskip (ii) We have $\Phi_p(x) = x^{p-1} + \dots + x + 1$, and $\Phi_{p^k}(x) = \underset{i = 0,\dots,p-1} \sum x^{i p^{k-1}}$ for
all integers $k \geqslant 1$ (see for instance \cite{L}, Chap. IV, \S 1, or \cite {FT}, Chap. VI, \S 1). Hence $\Phi_{p^k}(1) = p$
and $\Phi_{p^k}(-1) = 1$ for all $r \geqslant 1$. We have $\Phi_{2p^k}(x) = \Phi_{p^k}(-x)$, hence $\Phi_{2p^k}(1) = 1$ and
$\Phi_{2p^k}(-1) = p$
 for all $k \geqslant 1$. 
 
 \medskip (iii) Suppose that $m$ is divisible by at least two distinct prime numbers. Then $\Phi_m(1) = 1$ by \cite {W}, Proposition 2.8.
 We have $\Phi_m(-1) = \Phi_{2m}(1)$, therefore the same result implies that $\Phi_m(-1) = 1$.
 
 \medskip (iv) If $\Phi_m(1) = \Phi_m(-1) = 1$, then by (i) and (ii) the integer $m$ is not of the form $p^k$ or $2p^k$ for some prime number $p$.
 Therefore  $m$ is divisible by $4$ and by an odd prime number, or by two distinct odd prime numbers. This implies that $\varphi(m)$ is divisible by $4$, hence ${\rm deg}(\Phi_m) \equiv 0 \ {\rm (mod \ 4)}$.
 
 \begin{lemma}\label{divisible by 4}  Let $C = \Phi_m^r$ where 
$m, r$ are integers with $m \geqslant 3$ and $r \geqslant 1$. 
 If $C(1)$ and $C(-1)$ are both squares, then 
 ${\rm deg}(C)$ is divisible by 4.
 
 \end{lemma}
 
 \noindent
 {\bf Proof.} If $C(1) = C(-1) = 1$, then   Lemma \ref{1 and -1} (iv) implies that  ${\rm deg}(C)$ is divisible by 4.
Suppose that $C(1)C(-1) \not = 1$. Then by Lemma \ref{1 and -1} we have $m = p^k$ or $m = 2p^k$ for some prime number $p$, and
hence $\Phi_m(1) = p$ or $\Phi_m(-1) = p$ (see Lemma \ref{1 and -1} (i) and (ii)). This implies that $r$ is even, and since 
${\rm deg}(\Phi_m)$ is even, the degree of $C$ is divisible by 4.

 \bigskip
We now recall some notation from the introduction :

\medskip
$\bullet$ $C = \Phi_m^r$ where 
$m, r$ are integers with $m \geqslant 3$ and $r \geqslant 1$,

\medskip
$\bullet$ $R,S \geqslant 0$ are integers  such that  $R  \equiv S \ {\rm (mod \ 8)}$ and ${\rm deg}(C) < R + S$.
Set $N = R + S$.


\medskip

Theorem 2  is a consequence of lemmas \ref{P2 1} and \ref{P2 2}  :

\begin{lemma}\label{P2 1} Let $L$ be an even unimodular lattice of signature $(R,S)$, and let $t : L \to L$ be an isometry
with characteristic polynomial $C(x)(x-1)^{N - {\rm deg}(C)}$. Let $(c,{\rm deg}(C)-c)$ be the signature of the
sublattice ${\rm Ker}(C(t))$. Then

\medskip
{\rm (i)} $C(-1)$ is a square.

\medskip
{\rm (ii)} If $C(1) = 1$, then  ${\rm deg}(C)    \equiv \ 2c  \ {\rm (mod \ 8)}$.

 \medskip Moreover, if $C(1) =1$, then the sublattice ${\rm Ker}(C(t))$ is unimodular.

\end{lemma}

\medskip
 
 \noindent
 {\bf Proof.} Set $F(x) = C(x) (x-1)^{N-{\rm deg}(C)}.$ 
 Then by Lemma \ref{condition c1} the
 polynomial $F$ satisfies condition (C 1); this implies that $|F(-1)|$ is a square. Note that ${N-{\rm deg}(C)}$ is even, hence $|C(-1)|$ is a square. 
 Since $C$ is a power of a cyclotomic polynomial, we have $C(-1) \geqslant 0$ (cf. Lemma \ref{1 and -1}), hence $C(-1)$ is a square,
 and hence  (i) holds.
 
 \medskip
 Set $L_1 = {\rm Ker}(C(t))$, and let $L_2$ be the orthogonal complement of $L_1$
 in $L$.  If $C(1) = 1$, then the polynomials
 $x-1$ and $C$ are relatively prime over ${\bf Z}$ (i.e. the resultant of $x-1$ and $C$ is equal to 1); this implies that $L = L_1 \oplus L_2$, and hence the lattices $L_1$ and $L_2$
 are both even and unimodular. Therefore ${\rm deg}(C) - c    \equiv \ c  \ {\rm (mod \ 8)}$, hence
  ${\rm deg}(C)    \equiv \ 2c  \ {\rm (mod \ 8)}$, and therefore (ii) holds. 
  
   \begin{lemma}\label{P} 
{\rm (i)} If $C(1)$ and $C(-1)$ are both squares, then ${\rm deg}(C)    \equiv \ 0  \ {\rm (mod \ 4)}$. 

\medskip
{\rm (ii)} If $C(1)$ and $C(-1)$ are both squares, then Condition {\rm (C 1)} holds for $C$. 

\medskip
Let $c \geqslant 0$ be an even integer such that $c \leqslant {\rm deg}(C)$, $c \leqslant R$ and ${\rm deg}(C) - c \leqslant S$. We have

\medskip
{\rm (iii)} If $N = {\rm deg}(C) + 2$ and  ${\rm deg}(C)    \equiv \ 0  \ {\rm (mod \ 4)}$, then ${\rm deg}(C)    \equiv \ 2c  \ {\rm (mod \ 8)}$. 

\end{lemma}

\noindent
{\bf Proof.} If $C(1)$ and $C(-1)$ are both squares, then Lemma \ref{divisible by 4} implies 
that ${\rm deg}(C)$ is divisible by 4, hence (i) holds. Set ${\rm deg}(C) = 2n$; then $n$ is even, hence $C(1)$, $C(-1)$ and $(-1)^nC(1)$ and $C(-1)$
are all squares, therefore condition (C 1) holds for $C$, and this implies (ii).

\medskip
Let us prove (iii). Since ${\rm deg}(C)$ is divisible by 4 and ${\rm deg}(C) = N - 2$, we have
$N  \equiv \ 2  \ {\rm (mod \ 4)}$, and therefore $R$ and $S$ are both odd integers. Since ${\rm deg}(C) = N - 2$, 
the inequalities $c \leqslant R$ and ${\rm deg}(C) - c \leqslant S$ imply that $R-2 \leqslant c \leqslant R$; since $R$ and $S$ are odd, this implies that
$c = R -1$ and ${\rm deg}(C) - c = S-1$. W have $R  \equiv S \ {\rm (mod \ 8)}$ by hypothesis, hence 
${\rm deg}(C) - c   \equiv \ c  \ {\rm (mod \ 8)}$, as claimed.

  \begin{lemma}\label{P2 2} Let $c \geqslant 0$ be an even integer such that $c \leqslant {\rm deg}(C)$, $c \leqslant R$ and ${\rm deg}(C) - c \leqslant S$. Suppose that the following conditions hold

\medskip
{\rm (i)} $C(-1)$ is a square.

\medskip
{\rm (ii)} If $C(1) = 1$, then  ${\rm deg}(C)    \equiv \ 2c  \ {\rm (mod \ 8)}$. 

\medskip
Then there exists an even unimodular lattice $L$ 
 of signature $(R,S)$ and an isometry $t : L \to L$ such that
 
 \medskip
 $\bullet$ The characteristic polynomial of $t$ is $C(x)(x-1)^{N - {\rm deg}(C)}$,
 
 \medskip
 $\bullet$ The signature of the sublattice ${\rm Ker}(C(t))$ is $(c,{\rm deg}(C)-c)$,
 

  \end{lemma}
  
  \noindent
  {\bf Proof.} Set  $F(x) = C(x) (x-1)^{N-{\rm deg}(C)}.$
By (i), the polynomial $F$ satisfies condition (C 1). 

\medskip
Suppose that $C(1) > 1$, and note that by Lemma \ref{1 and -1} this implies that
$m = p^k$ for some prime number $p$. If $C(1) > 1$ and ${\rm deg}(C) < N - 2$, then with the notation of \S \ref{isometries} we have $n^+ \not = 2$,
hence $\Pi_{\Phi_m(x),x-1} = \{p\}$. Suppose now that $C(1) > 1$, that ${\rm deg}(C) = N - 2$, and that
$C(1)$ is not a square. Then  we have $D_0 \not = -1$ in ${\bf Q}_p^{\times}/{\bf Q}_p^{\times 2}$, hence 
$\Pi_{\Phi_m(x),x-1} = \{p\}$ in this case as well. Therefore in both cases we have $G_F = 0$. 
By Theorem \ref{signature map theorem}, there exists an even, unimodular lattice $L$ of signature $(R,S)$ and an isometry
$t : L \to L$ with characteristic polynomial $F$ and signature map $\tau$ satisfying $\tau(C) = (c,{\rm deg}(C) - c)$. Let $L_1 = {\rm Ker}(C(t))$, and let $L_2$ be the
sublattice of $L$ of fixed points by $t$. 

\medskip
It remains to consider the cases where $C(1)$ is a square, and either $C(1) = 1$, or ${\rm deg}(C) = N - 2$.

\medskip
If $C(1)$ is a square, then by Lemma \ref{P} condition (C 1) holds for $C$. Moreover, we have ${\rm deg}(C) - c   \equiv \ c  \ {\rm (mod \ 8)}$. 
If $C(1) = 1$, this follows from (ii), and if ${\rm deg}(C) = N - 2$, from Lemma \ref{P} (i) and (iii).

\medskip
The polynomial $C$ is  a power of an irreducible polynomial, hence the group $G_C$ is trivial, and therefore in both cases we can apply  Theorem \ref{basic theorem}, and conclude that there exists
an even unimodular lattice $L_1$ of signature $(c,{\rm deg}(C) - c)$ and an isometry $t_1 : L_1 \to L_1$ of characteristic polynomial $C$ (note that this also
follows from \cite{BT}, Theorem A). Let $L_2$ be an even unimodular lattice of signature $(R-c,S-{\rm deg}(C)-c)$,  and let $t_2 : L_2 \to L_2$ be the identity. Set $L = L_1 \oplus L_2$, and $t = (t_1,t_2)$.

\medskip

The lattice $L$ is even unimodular and of signature $(R,S)$, and $t$ is an isometry of
$L$ with the required properties. This completes the proof of the lemma.

\bigskip
The following is a reformulation of Theorem 2 of the introduction : 
 
 \begin{theo}\label{Theorem 2}
  There exists an even unimodular lattice $L$ 
 of signature $(R,S)$ and an isometry $t : L \to L$ such that
 
 \medskip
 $\bullet$ The characteristic polynomial of $t$ is $C(x)(x-1)^{N - {\rm deg}(C)}$,
 
 \medskip
 $\bullet$ The signature of the sublattice ${\rm Ker}(C(t))$ is $(c,{\rm deg}(C)-c)$,
 
 \bigskip if and only if the following conditions are satisfied : 

\bigskip
{\rm (i)} $C(-1)$ is a square.

\medskip
{\rm (ii)} If $C(1) = 1$, then  ${\rm deg}(C)    \equiv \ 2c  \ {\rm (mod \ 8)}$.

 \end{theo}
 
 \medskip
 
 \noindent
 {\bf Proof.} This is an immediate consequence of Lemmas \ref{P2 1} and \ref{P2 2}. 
 
 \begin{notation} If $L$ is a lattice, we denote by $L^{\sharp}$ its dual lattice, and set $\Delta(L) = L^{\sharp}/L$.
 
 \end{notation}
 
 \begin{defn}  Let $p$ be
 a prime number. We say that a lattice $L$ is {\it $p$-elementary} if $p \Delta(L)$ = 0. 
 
 \end{defn}
 
 \begin{prop}\label{elementary} Let $L$ be a unimodular lattice of rank $N$ and let $t : L \to L$ be an isometry with
 characteristic polynomial $C(x)(x-1)^{N - {\rm deg}(C)}$. Set $L_C = {\rm Ker}(C(t))$, and let $L_0$ be the
 orthogonal complement of $L_C$ in $L$. Then we have
 
 \medskip
{\rm (i)} If $C(1) = 1$, then $L_C$ and $L_0$ are both unimodular.
 
 \medskip
 {\rm (ii)} If $C(1) > 1$, then $L_C$ and $L_0$ are both $p$-elementary, where $p$ is such that $\Phi_m(1) = p$. 
 
 \end{prop}
 
 \noindent
 {\bf Proof.} (i) Lemma \ref{P2 1} implies that $L_C$ is unimodular, hence $L_0$ is also unimodular.
 
 \medskip (ii) The action of $t$ endows $L$ with a structure of ${\bf Z}[\Gamma]$-module with $\Gamma$ infinite cyclic; this
 action stabilizes $L_C$ and $L_0$, hence also $\Delta(L_C)$ and $\Delta(L_0)$. Since $L$ is unimodular, the ${\bf Z}[\Gamma]$-modules
 $\Delta(L_C)$ and $\Delta(L_0)$ are isomorphic. 
 
 \medskip
 Lemma \ref{1 and -1} implies that $C = \Phi_m^r$ with $m = p^k$ for some integer $k$. Therefore $t$ acts on $L_C$ by
 $t(x) = \alpha.x$ where $\alpha$ is the image of $x$ in ${\bf Z}[x]/(\Phi_{p^k})$. On the other hand, $t$ acts as the identity on $L_0$.
 Since the ${\bf Z}[\Gamma]$-modules
 $\Delta(L_C)$ and $\Delta(L_0)$ are isomorphic, this implies that $p \Delta(L_C) = 0$ and $p \Delta(L_0) = 0$.

 \bigskip

\section{Automorphisms acting trivially on Picard lattices}\label{Vorontsov}

We now prove Proposition 2 of the introduction. Let $m, r$ be integers with $m \geqslant 3$ and $r \geqslant 1$, and let $C = \Phi_m^r$. Assume
that ${\rm deg}(C) \leqslant 20$. 

\begin{prop}\label{Kondo} There exists an automorphism $a : X \to X$ of a  projective $K3$ surface $X$ such that $a^*$ is the identity on $S_X$ and that the characteristic polynomial of
the restriction of $a^*$ to $T_X$ is equal to $C$ if and only if the following conditions hold

\medskip
{\rm (i)} $C(-1)$ is a square.

\medskip
{\rm (ii)} If $C(1) = 1$, then  ${\rm deg}(C)  \equiv \ 4  \ {\rm (mod \ 8)}$. 

\medskip If $C(1) = 1$, then the lattice $T_X$ is unimodular. 
Moreover, the $K3$ surface is unique up to isomorphism if and only if $C$ is a cyclotomic polynomial {\rm (i.e. $r = 1$)}. 

\end{prop}

\noindent
{\bf Proof.} Let $a : X \to X$  be an automorphism of a projective $K3$ surface such that $a^*|S_X$ is the identity, and that the
characteristic polynomial of $a^*|T_X$ is equal to $C$. Applying Lemma \ref{P2 1} with $L = H^2(X,{\bf Z})$, $t = a^*$, $N = 22$, $R = 3$, $S = 19$ and $c = 2$
shows that (i) and (ii) hold.

\medskip
Conversely, assume that (i) and (ii) hold,
and set $$F(x) = C(x) (x-1)^{22-{\rm deg}(C)}.$$ By Lemma \ref{P2 2} with $N = 22$, $R = 3$, $S = 19$ and $c = 2$ there
exists an 
even unimodular lattice $L$ of signature $(3,19)$ and an isometry
$t : L \to L$ with characteristic polynomial $F$ such that the signature of  the sublattice $L_1 = {\rm Ker}(C(t))$ is $(2,{\rm deg}(C) - 2)$. Let $L_1 = {\rm Ker}(C(t))$, and let $L_2$ be the
sublattice of $L$ of fixed points by $t$.

\medskip Let $V \subset L_1 \otimes_{\bf Z} {\bf R}$ be a $2$-dimensional subspace of signature $(2,0)$ and  stable by $t$; for a generic choice of $V$, the 
intersection of $L$ with the orthogonal of $V$ is equal to $L_2$. The restriction of $t$ to $V$ has determinant 1. Since $t_2$ is the identity, it is
positive in the sense of McMullen, \cite{Mc3}, \S 2. By \cite{Mc3}, Theorem 6.1, there exists a projective $K3$ surface $X$ and an automorphism $a : X \to X$
such that $a^* = t$, that $S_X = L_2$ and $T_X = L_1$. Therefore $a^*$ is the identity on $S_X$ and the restriction of $a^*$ to $T_X$ has characteristic polynomial $C$. If moreover $C(1) = 1$, then by Lemma \ref{P2 1} the lattice $T_X$ is unimodular. 

\medskip If $r = 1$, then the
uniqueness of the $K3$ surface up to isomorphism follows from a result of Brandhorst, \cite{Bra}, Theorem 1.2. If $r > 1$, then varying the choice of the subspace $V$ gives
rise to an infinite family of $K3$ surfaces. 
This completes the proof of the proposition.

\begin{coro}\label{p-primary} Let $a : X \to X$ be an automorphism of a  projective $K3$ surface $X$ such that $a^*$ is the identity on $S_X$ and that the characteristic polynomial of
the restriction of $a^*$ to $T_X$ is equal to $C$. 

\medskip Then one of the following holds :

\medskip
{\rm (i)} $T_X$ and $S_X$ are unimodular.

\medskip
{\rm (ii)} $T_X$ and $S_X$ are $p$-elementary, where $p$ is a prime number such that $\Phi_m(1) = p$.

\end{coro}

\noindent
{\bf Proof.} This is an immediate consequence of Proposition \ref{elementary}. 

\medskip

Note that this implies that $p \leqslant 19$, hence we have the following :

\begin{coro}\label{19} If a  projective $K3$ surface $X$ has a non-trivial automorphism that induces the identity on $S_X$, then the lattices $T_X$ and $S_X$ are either both unimodular, or both $p$-elementary with 
$p \leqslant 19$. 

\end{coro}

\section{Cyclotomic polynomials and isometries of finite order}\label {P 1 section}

We keep the notation of \S \ref{P 2}; in particular, 
$C = \Phi_m^r$ where 
$m, r$ are integers with $m \geqslant 3$ and $r \geqslant 1$. 
The 
following result implies Theorem 1 of the introduction; it is actually a {\it strengthening} of Theorem 1, since 
it implies the existence of an
isometry is of order $m$. 

\begin{theo}\label{theorem 1} Suppose that $C(1)C(-1) > 1$. Then there exists an even, unimodular lattice 
$L$ 
 of signature $(R,S)$ and an isometry $t : L \to L$ of order $m$  such that
 
 \medskip
 $\bullet$ The characteristic polynomial of $t$ is divisible by $C$,
 
 \medskip
 $\bullet$ The signature of the sublattice ${\rm Ker}(C(t))$ is $(c,{\rm deg}(C)-c)$.

\end{theo}

\noindent
{\bf Proof.} If $C(1) > 1$ and $C(-1)$ is a square, then this follows from Theorem \ref{Theorem 2}. 

\medskip
Suppose that $C(-1) > 1$, and 
 that $C(1)$ is a square. In this case, set $F(X) = C(X)(X+1)^{N - {\rm deg}(C}$. Since $C(1)$ is a square, the polynomial
$F$ satisfies Condition (C 1). 

\medskip If 
${\rm deg}(C) < N - 2$, then with the notation of \S \ref{isometries} we have $n^- \not = 2$,
hence $\Pi_{\Phi_m(x),x-1} \not = \varnothing$; this implies that $G_F = 0$. Suppose now that ${\rm deg}(C) = N - 2$, and that
$C(-1)$ is not a square.  Note that by Lemma \ref{1 and -1} this implies that $m = 2p^k$ for some odd prime number $p$; with the notation of \S \ref{isometries}, we have $D_0 \not = -1$ in ${\bf Q}_p^{\times}/{\bf Q}_p^{\times 2}$, hence 
$\Pi_{\Phi_m(x),x-1} = \{p\}$. This implies that $G_F = 0$ in this case as well.
By Theorem \ref{signature map theorem}, there exists an even, unimodular lattice $L$ of signature $(R,S)$ and an isometry
$t : L \to L$ with characteristic polynomial $F$ and signature map $\tau$ satisfying $\tau(C) = (c,{\rm deg}(C) - c)$. 

\medskip Suppose now that $C(-1)$ is a square and  ${\rm deg}(C) = N - 2$; then Lemma \ref{P} implies that condition (C 1) holds for $C$ and 
that ${\rm deg}(C) - c   \equiv \ c  \ {\rm (mod \ 8)}$.
The polynomial $C$ is  a power of an irreducible polynomial, hence the group $G_C$ is trivial, and therefore  we can apply  Theorem \ref{basic theorem}, and conclude that there exists
an even unimodular lattice $L_1$ of signature $(c,{\rm deg}(C) - c)$ and an isometry $t_1 : L_1 \to L_1$ of characteristic polynomial $C$ (note that this also
follows from \cite{BT}, Theorem A). Let $L_2$ be an even unimodular lattice of signature $(R-c,S-{\rm deg}(C)-c)$,  and let $t_2 : L_2 \to L_2$ be the identity. Set $L = L_1 \oplus L_2$, and $t = (t_1,t_2)$; then $t : L \to L$ has the required properties.

\medskip 

Finally, suppose that $C(1)$ and $C(-1)$ are both non-squares. In this case, $C(1) = C(-1) = 2$, and $C = \Phi_{2^k}$ for some integer $k$ (see
Lemma \ref{1 and -1}). Suppose first that $N > {\rm deg}(C) + 2$. Set $F(x) = C(x) (x+1)^2 (x-1)^{N - {\rm deg}(C) - 2}$. We have
$\Pi_{C,x-1} = \Pi_{x-1,x+1} = \{2\}$ (see \cite{B 22}, \S 7 and \S 12), hence $G_F = 0$. Therefore by \cite {B 22}, Corollary 12.3 there
exists exists an even, unimodular lattice $L$ of signature $(R,S)$ and an isometry
$t : L \to L$ with characteristic polynomial $F$ and signature map $\tau$ satisfying $\tau(C) = (c,{\rm deg}(C) - c)$. Assume now that
$N = {\rm deg}(C) + 2$, and set $F(x) = C(x) (x+1) (x-1)$. By Takada \cite{T}, Theorem 6.11, there exist a lattice $L$ and an isometry $t : L \to L$
with the required properties. This concludes the proof of the theorem. 

\medskip
It remains to treat the case where $C(1) = C(-1) = 1$;  if ${\rm deg}(C)    \equiv \ 2c  \ {\rm (mod \ 8)}$, then
Theorem \ref{Theorem 2} implies the following

\begin{coro}\label{coro 2}
Suppose that $C(1) = 1$ and that   ${\rm deg}(C)    \equiv \ 2c  \ {\rm (mod \ 8)}$. 
Then
there exists an even unimodular lattice $L$ 
 of signature $(R,S)$ and an isometry $t : L \to L$ of order $m$ such that
 
 \medskip
 $\bullet$ The characteristic polynomial of $t$ is divisible by $C$,
 
 \medskip
 $\bullet$ The signature of the sublattice ${\rm Ker}(C(t))$ is $(c,{\rm deg}(C)-c)$.

 \end{coro}

\medskip
The condition $C(1) = C(-1) = 1$ implies that ${\rm deg}(C)    \equiv \ 0  \ {\rm (mod \ 4)}$ (see Lemma \ref{1 and -1} (i)),
and $c$ is an even integer; hence we have either ${\rm deg}(C)    \equiv \ 2c  \ {\rm (mod \ 8)}$, or
${\rm deg}(C)    \equiv \ c  \ {\rm (mod \ 8)}$. The first case is covered by Corollary \ref{coro 2}, therefore we have
the following two cases to consider

\medskip
(a)  $c    \equiv \ 0 \ {\rm (mod \ 4)}$ and ${\rm deg}(C)    \equiv \ 4  \ {\rm (mod \ 8)}$

\medskip
(b) $c    \equiv \ 2 \ {\rm (mod \ 4)}$ and ${\rm deg}(C)    \equiv \ 0  \ {\rm (mod \ 8)}$.

\medskip We treat these cases in the next sections; the following results will be useful~:

\begin{lemma}\label{indefinite} Let $t : L \to L$ be an isometry of a lattice $L$ such that 

 \medskip
 $\bullet$ The characteristic polynomial of $t$ is divisible by $C$,
 
 \medskip
 $\bullet$ The signature of the sublattice ${\rm Ker}(C(t))$ is $(c,{\rm deg}(C)-c)$.
 
 \medskip
 If ${\rm deg}(C)  \equiv \ 0 \ {\rm (mod \ 4)}$ and $c    \equiv \ 2 \ {\rm (mod \ 4)}$, then $L$ is indefinite.

\end{lemma}

\noindent
{\bf Proof.} Indeed, the sublattice ${\rm Ker}(C(t))$ is indefinite : since ${\rm deg}(C)  \equiv \ 0 \ {\rm (mod \ 4)}$  and
$c    \equiv \ 2 \ {\rm (mod \ 4)}$, we have $c \not = {\rm deg}(C)$ and $c \not = 0$.

 \begin{prop}\label{pi} Let $m \geqslant 3$ be an integer and let $p$ be a prime number that does not divide $m$. The following
 are equivalent
 
 \medskip
 {\rm (a)} The polynomial $\Phi_m$ has a symmetric irreducible factor mod $p$.
 
 \medskip
 {\rm (b)} The prime ideals of ${\bf Q}(\zeta_m + \zeta_m^{-1})$ above $p$ are inert in ${\bf Q}(\zeta_m)$.
 
 \medskip
 {\rm (c)} The subgroup of $({\bf Z}/m {\bf Z})^{\times}$ generated by $p$ contains $-1$. 
 
 \end{prop}
 
 \noindent
 {\bf Proof} The equivalence of (a) and (b) follows from \cite{W}, Proposition 2.14. Let us prove that (b) and (c) are equivalent. 
 Let $G$ be the Galois group  of ${\bf Q}(\zeta_m)/{\bf Q}$, and let  $P$ be a prime ideal of ${\bf Q}(\zeta_m)$ above $p$. The decomposition group $G_P$ is by definition $\{ g \in G \ | \ g(P) = P \}$. 
 Since $G$ is abelian, this group only depends on the prime number $p$; set $G_P = G_p$. Condition (b) holds if and only
 if  the element of $G$ induced by $\zeta_m \to \zeta_m^{-1}$ is contained in $G_p$. Let $f : G \to ({\bf Z}/m {\bf Z})^{\times}$ be an
 isomorphism; then $f(G_p)$ is the subgroup of $({\bf Z}/m {\bf Z})^{\times}$ generated by $p$. This
 implies the equivalence of (b) and (c).

 \begin{coro}\label {utile} Let $m, p$ be distinct prime numbers. If $p$ is not a square modulo $m$ then the polynomial $\Phi_m$ has a symmetric irreducible factor mod $p$.
 
 \end{coro}
 
 \noindent
 {\bf Proof.} We have $p^{{m-1}\over 2} = \pm 1$. If $p^{{m-1}\over 2} = 1$, then $p$ is a square modulo $m$, hence
 $p^{{m-1}\over 2} = - 1$. This implies that the subgroup of $({\bf Z}/m {\bf Z})^{\times}$ generated by $p$ contains $-1$,
 and hence by Proposition \ref{pi} the  polynomial $\Phi_m$ has a symmetric irreducible factor mod $p$.
 

 
 


\medskip
We start by noting that if $N$ is sufficiently large, then Property (P 1') holds.

\begin{prop}\label{large N} Suppose that $C(1) = C(-1)  = 1$. Let $p$ be a prime number such that $\Pi_{\Phi_{mp},\Phi_m} = \{p\}$. If $N > {\rm deg}(C)+  \varphi(mp) $,
then there exists an even unimodular lattice $L$ 
 of signature $(R,S)$ and an isometry $t : L \to L$ of order $mp$ such that
 
 \medskip
 $\bullet$ The characteristic polynomial of $t$ is divisible by $C$,
 
 \medskip
 $\bullet$ The signature of the sublattice ${\rm Ker}(C(t))$ is $(c,{\rm deg}(C)-c)$.

\end{prop}

\noindent
{\bf Proof.} Set $F(x) = C(x)\Phi_{mp}(x)(x-1)^k$, with $k = N - {\rm deg}(C) - \varphi(mp) $. The polynomial $F$ satisfies
condition (C 1). Since $\Pi_{\Phi_{mp},\Phi_m} = \{p\}$, we have $G_F = 0$. Therefore by Theorem \ref{signature map theorem}
there exists a lattice $L$ and an isometry $t$ with the required properties.

\medskip

Note that Proposition \ref{pi}  implies that there exist infinitely many prime numbers $p$ such that $\Pi_{\Phi_{mp},\Phi_m} = \{p\}$.
In the following sections, we give conditions on $N$ for
the existence of an isometry {\it of order $m$}.

\section{$C(1) = C(-1) = 1$ and $c    \equiv \ 0 \ {\rm (mod \ 4)}$}\label{0 section}

We keep the notation of the previous sections, and we assume in addition that $C(1) = C(-1) = 1$; this implies that 
${\rm deg}(C)    \equiv \ 0  \ {\rm (mod \ 4)}$ (see Lemma \ref{1 and -1} (i)). Suppose that  $c    \equiv \ 0 \ {\rm (mod \ 4)}$.

\medskip
Corollary \ref{coro 2} implies that if ${\rm deg}(C)    \equiv \ 0  \ {\rm (mod \ 8)}$, then there exists an even, unimodular lattice $L$ 
of signature $(R,S)$ having an isometry $t : L \to L$ of order $m$ such that the signature of ${\rm Ker}(C(t))$ is $(c, {\rm deg}(C) - c)$.

\medskip
Suppose that  ${\rm deg}(C)    \equiv \ 4  \ {\rm (mod \ 8)}$; then $r$ is odd, and Lemma \ref{1 and -1} implies that $m$ is of one of the following forms

\medskip
$\bullet$ $m = 4p^k$ where $p$ is a prime number with $p   \equiv \  3 \  {\rm (mod \ 4)}$ and $k \geqslant 1$ is an integer;

\medskip
$\bullet$  $m = p^kq^s$ where $p$ and $q$ are distinct prime numbers with $\equiv \ 3 \ {\rm (mod \ 4)}$ and $k,s \geqslant 1$ are integers;

\medskip
$\bullet$  $m = 2p^kq^s$ where $p$ and $q$ are distinct prime numbers with $\equiv \ 3 \ {\rm (mod \ 4)}$ and $k,s \geqslant 1$ are integers.

\begin{lemma}\label{N big}  $N \geqslant  {\rm deg}(C) + 4$.

\end{lemma}

\noindent
{\bf Proof.} Let us show that $N \not = {\rm deg}(C) + 2$. Set $c' = {\rm deg}(C) - c$. We have $c \leqslant R$, $c' \leqslant S$ and
$N = R + S$, ${\rm deg}(C) = c + c'$; moreover, $c$ and $c'$ are even. Therefore if $N = {\rm deg}(C) + 2$, then $R = c + 1$
and $S = c' + 1$. We have  $R    \equiv \ S \ {\rm (mod \ 8)}$, hence this implies that $c    \equiv \ c' \ {\rm (mod \ 8)}$; but
${\rm deg}(C) = c + c'$ is congruent to $4$ {\rm (mod \ 8)}, so this is impossible. Since $N$ and $ {\rm deg}(C)$ are both even,
this implies that $N \geqslant  {\rm deg}(C) + 4$, as claimed.

\begin{prop}\label{4p} Suppose that $m = 4p^k$ where $p$ is a prime number with $p   \equiv \  3 \  {\rm (mod \ 4)}$ and $k \geqslant 1$ is an integer.
Then there exists an even unimodular lattice $L$ 
 of signature $(R,S)$ and an isometry $t : L \to L$ of order $m$ such that
 
 \medskip
 $\bullet$ The characteristic polynomial of $t$ is divisible by $C$,
 
 \medskip
 $\bullet$ The signature of the sublattice ${\rm Ker}(C(t))$ is $(c,{\rm deg}(C)-c)$.

\end{prop}

\noindent
{\bf Proof.}  Set $F = C \Phi_4^2(x-1)^{N-{\rm deg}(C) -4}$; note that  Lemma \ref{N big} implies that $N-{\rm deg}(C) -4 \geqslant 0$, and that $F$ satisfies condition (C 1). 
Since $p   \equiv \  3 \  {\rm (mod \ 4)}$, 
$\Phi_4$ is irreducible {\rm mod  p}, and therefore $\Pi_{\Phi_m,\Phi_4} = \{p\}$. This implies that $G_F = 0$, and hence
by  Theorem \ref{signature map theorem}
there exists a lattice $L$ and an isometry $t$ with the required properties. 

\medskip If $p$ and $q$ are distinct prime numbers $\equiv \ 3 \ {\rm (mod \ 4)}$, then by quadratic reciprocity either 
$p$ is a square modulo $q$, or $q$ is a square modulo $p$, and these cases are mutually exclusive. Therefore 
we may assume that $p$ is a square modulo $q$.

\begin{prop}\label{p and 2p} Suppose that $m = p^kq^s$ or $2p^kq^s$ where $p$ and $q$ are distinct prime numbers with $p, q \equiv \ 3 \ {\rm (mod \ 4)}$ and $k,s \geqslant 1$ are integers, and assume that $p$ is a square modulo $q$. Suppose that $N \geqslant {\rm deg}(C) + (p-1)p^{k-1} + 2$.

\medskip

Then there exists an even unimodular lattice $L$ 
 of signature $(R,S)$ and an isometry $t : L \to L$ of order $m$ such that
 
 \medskip
 $\bullet$ The characteristic polynomial of $t$ is divisible by $C$,
 
 \medskip
 $\bullet$ The signature of the sublattice ${\rm Ker}(C(t))$ is $(c,{\rm deg}(C)-c)$,
 

\end{prop}

\noindent
{\bf Proof.} Set $N' = N - {\rm deg}(C) - (p-1)p^{k-1} - 2$; set $F(x) = C(x)\Phi_{p^k}(x)(x-1)^{N'}$
if $m =  p^kq^s$, and $F(x) = C(x)\Phi_{2p^k}(x)(x+1)^{N'}$ if $m = 2p^kq^s$. 
The polynomial $F$ satisfies condition (C 1), and Lemma \ref{1 and -1}
implies that $G_F = 0$. Theorem \ref{signature map theorem} implies that
there exists a lattice $L$ and an isometry $t$ with the required properties. 




\section{$C(1) = C(-1) = 1$ and $c    \equiv \ 2 \ {\rm (mod \ 4)}$}\label{2 section}

We keep the notation of \S \ref{P 2}; in particular, 
$C = \Phi_m^r$ where 
$m, r$ are integers with $m \geqslant 3$ and $r \geqslant 1$. The case where $C(1)C(-1) > 1$ is covered by Theorem \ref{theorem 1}. 

\medskip
Assume now that $C(1) = C(-1) = 1$, and that
$c    \equiv \ 2 \ {\rm (mod \ 4)}$. Since $C(1) = C(-1) = 1$, by Lemma \ref{1 and -1} we have ${\rm deg}(C)    \equiv \ 0  \ {\rm (mod \ 4)}$; hence
Lemma \ref{indefinite} implies that if $t : L \to L$ is an isometry of a lattice such that he signature of ${\rm Ker}(C(t))$ is $(c, {\rm deg}(C) - c)$,
then $L$ is indefinite. This implies that $R > 0$ and $S > 0$; recall that 
since $R  \equiv S \ {\rm (mod \ 8)}$,  there exists up to isomorphism a unique even, unimodular lattice of signature $(R,S)$ (see for
instance \cite {S}, Chap. V); we denote it by $\Lambda_{R,S}$. 

\medskip
In the applications to $K3$ surfaces, we have $R = 3$, $S = 19$ and $c = 2$.





\medskip
Note that the case where ${\rm deg}(C)    \equiv \ 4 \ {\rm (mod \ 8)}$ was already handled in 
Corollary \ref{coro 2} :

\begin{prop} If  ${\rm deg}(C)    \equiv \ 4 \ {\rm (mod \ 8)}$ then the lattice $\Lambda_{R,S}$ has an isometry $t$  of order $m$
such that the signature of ${\rm Ker}(C(t))$ is $(c, {\rm deg}(C) - c)$.

\end{prop} 

\noindent
{\bf Proof.} Indeed, we are assuming that $c    \equiv \ 2 \ {\rm (mod \ 4)}$, hence the hypothesis 
${\rm deg}(C)    \equiv \ 4 \ {\rm (mod \ 8)}$  implies that ${\rm deg}(C)    \equiv \ 2c \ {\rm (mod \ 8)}$. Therefore by Corollary \ref{coro 2}
there exists an even, unimodular lattice $L$ and an isometry $t : L \to L$ of order $m$ such that the
signature of ${\rm Ker}(C(t))$ is $(c, {\rm deg}(C) - c)$. By Lemma \ref{indefinite}, such a lattice is indefinite, hence $L$ is
isomorphic to $\Lambda_{R,S}$. This concludes the proof of the lemma.

\bigskip
Suppose that  ${\rm deg}(C)    \equiv \ 0  \ {\rm (mod \ 8)}$. Recall that $N = R + S$. Using the results of \cite{B 22} (in particular, 
Theorem \ref{signature map theorem}) and Proposition \ref{pi} it is possible to determine the values of $N$ for which $\Lambda_{R,S}$
has an isometry $t$ of order $m$ such that the signature of ${\rm Ker}(C(t))$ is $(c, {\rm deg}(C) - c)$; since this would be
quite long, we only give some partial results that will be useful for the for the proof of Proposition 1.

\begin{prop}\label{15} Let $m = 2^n p$ with $n \geqslant 2$ and $p$ a prime number such that $p    \equiv \ 3,5 \ {\rm (mod \ 8)}$, or $m = pq$ with $p$ and $q$ distinct
prime numbers such that $p$ is not a square modulo $q$.  Suppose that $N \geqslant {\rm deg}(C) + p + 1$; set $M = N - {\rm deg}(C) - p + 1$
and $F(x) = C(x)\Phi_p(x)(x-1)^M$.

\medskip
Then $\Lambda_{R,S}$ has an isometry $t$ with characteristic polynomial $F$ such that the signature of ${\rm Ker}(C(t))$ is $(c, {\rm deg}(C) - c)$
and that the signature of ${\rm Ker}((t-1)^M)$ is $(1, M-1)$.

\end{prop}

\noindent
{\bf Proof.} Since $C(1) = C(-1) = 1$, the polynomial $F$ satisfies condition (C~1). We have $G_F = 0$; indeed,
by Corollary \ref{utile} we have $\Pi_{\Phi_m,\Phi_p} = \{2\}$ if $m = 2^n p$ and $\Pi_{\Phi_m,\Phi_p} = \{q\}$ if $m = pq$; moreover,
$\Pi_{\Phi_p(x),x-1} = \{p\}$. Therefore by Theorem \ref{signature map theorem} there exists an isometry with the required properties.

\begin{prop}\label{30} Let  $m = 2pq$ with $p$ and $q$ distinct
prime numbers such that $p$ is not a square modulo $q$.  Suppose that $N \geqslant {\rm deg}(C) + p + 1$; set $M = N - {\rm deg}(C) - p - 1$
and $F(x) = C(x)\Phi_{2p}(x)(x+1)^2(x-1)^M$.

\medskip
Then $\Lambda_{R,S}$ has an isometry $t$ with characteristic polynomial $F$ such that the signature of ${\rm Ker}(C(t))$ is $(c, {\rm deg}(C) - c)$
and that the signature of ${\rm Ker}((t-1)^M)$ is $(1, M-1)$.

\end{prop}

\noindent
{\bf Proof.} We have $\Pi_{\Phi_m,\Phi_{2p}} = \{q\}$ by Corollary \ref{utile}, and $\Pi_{\Phi_{2p}(x),x +1} = \{p\}$, $\Pi_{x+1,x-1} = \{2\}$,
hence $G_F = 0$.  The polynomial $F$ satisfies condition (C~1). Therefore by Theorem \ref{signature map theorem} there exists an isometry with the required properties.



 \section{Salem polynomials and isometries of lattices}\label{Salem section}
 
 A {\it Salem polynomial} is a monic irreducible polynomial $S \in {\bf Z}[X]$ such that 
$S(X) = X^{{\rm deg}(S)} S(X^{-1})$ and that $S$ has exactly two roots outside the unit circle, both positive
real numbers.

 \begin{example}\label{Sa}  Let $n$ be an integer $\geqslant 0$, and set
$$S_n(X) = X^6 - nX^5 - X^4 + (2n-1) X^3 - X^2 - nX +1.$$  The polynomials $S_n$ are Salem polynomials 
(see  \cite{Mc1}, \S 4, or \cite{GM}, \S 7, Example 1); we have $S_n(1) = -1$ and $S_n(-1) = 1$. 

\end{example} 

\medskip If $a : X \to X$ is an automorphism of a projective $K3$ surface, then the characteristic polynomial of $a^* : H^2(X,{\bf C}) \to H^2(X,{\bf C})$ 
is either a product of cyclotomic polynomials, or it is of the form $SC$, where $S$ is a Salem polynomial and $C$ a product of cyclotomic polynomials (see 
\cite{Mc1}, Theorem 3.2 and Corollary 3.3). 


\medskip

We recall some notions and results from \cite{B 22}, \S 7 and \S 12. 

\begin{notation}\label{notation} Let $S$ be a Salem polynomial such that $S(1) = -1$ and $S(-1) = 1$, and let $C$ be a cyclotomic polynomial. Let $\Pi_{S,C}$ be the set of prime numbers
$p$ such that $S \ {\rm (mod \ {\it p})}$ and $C  \ {\rm (mod \ {\it p})}$
 have a common irreducible symmetric 
factor in $ {\bf F}_p[x]$.

\end{notation}

\begin{example}\label{example} Let $S_2(x) = x^6 - 2x^5 - x^4 + 3x^3 - x^2 - 2x + 1$ (cf example \ref{Sa}) and 
let $C = \Phi_{60}$. The polynomials $S_2 \ {\rm (mod \ {359})}$ and $C \ {\rm (mod \ {359})}$
have the common irreducible factor $x^2 - 15x + 1$ 
 in $ {\bf F}_{359}[x]$; this polynomial is symmetric, hence $359 \in \Pi_{S,C}$.

\end{example} 

\begin{notation}\label{group} Let $F = SC$ for $S$ and $C$ as in notation \ref{notation}. We define a group $G_F$ as in \cite{B 22}, \S 7 (see also \S \ref{isometries}); we have $G_F = 0$ if $\Pi_{S,C} \not = \varnothing$, and $G_F = {\bf Z}/2{\bf Z}$ if $\Pi_{S,C} = \varnothing$. 

\end{notation}

\begin{prop}\label{proposition} Let $S$ be a Salem polynomial such that $S(1) = -1$ and $S(-1) = 1$, and let $C$ be a cyclotomic polynomial; set $F = SC$. Suppose that ${\rm deg}(F) = 22$,
that condition {\rm (C 1)} holds for $F$, and that $G_F = 0$. Then the lattice $\Lambda_{3,19}$ has an isometry $t$ of signature map $\tau$ satisfying $\tau(C) = (2,{\rm deg}(C) -2)$.

\end{prop} 

\noindent 
{\bf Proof.} This is a consequence of \cite{B 22}, Corollary 12.3. 

\begin{example}\label{for 60 cyclotomic}  Let $S(x) = S_2(x) = x^6 - 2x^5 - x^4 + 3x^3 - x^2 - 2x + 1$, and
$C = \Phi_{60}$; set $F = SC$. We have  $359 \in \Pi_{S,C}$ (cf. Example \ref{example}), therefore $G_F = 0$. Condition (C 1) holds for $F$, hence by
Proposition \ref{proposition} the lattice $\Lambda_{3,19}$ has an isometry $t$ with characteristic polynomial $F$ such that the signature of ${\rm Ker}(C(t))$ is $(2,14)$.

\end{example} 



\begin{example}\label{30} Let $C = \Phi_{30}^2$ and 
$S_1(x) = x^6 - x^5 - x^4 + x^3 - x^2 - x + 1$, as in Example \ref{Sa}. We have
$3 \in \Pi_{S_1,\Phi_{30}}$, hence $G_{S_1C}= 0$. Proposition \ref{proposition} implies that  $\Lambda_{3,19}$ has an isometry $t$ such that the signature of ${\rm Ker}(C(t))$ is $(2,14)$.

\end{example}

\begin{example}\label{10} Let $C = \Phi_{10}^4$. We have $S_0(x) = x^6 - x^4 - x^3 - x^2 + 1$ (cf Example \ref{Sa}) and
$3 \in \Pi_{S_0,\Phi_{10}}$, hence $G_{S_0C}= 0$. Proposition \ref{proposition} implies that  $\Lambda_{3,19}$ has an isometry $t$ such that the signature of ${\rm Ker}(C(t))$ is $(2,14)$.

\end{example}

\begin{remark} The sets $\Pi_{S,C}$ of the above examples were computed by PARI GT. 

\end{remark} 

\begin{notation} Let $C$ be a cyclotomic polynomial, and let $S_n$ be as in example \ref{Sa}. Let $N(C)$ be the set of integers $n \geqslant 0$
such that $\Pi_{C,S_n} \not = \varnothing$.

\end{notation}

\begin{example} Let $C = \Phi_{60}$. We have $0, 2, 5, 6, 7,... \in N(C)$. 

\end{example}

\begin{question} Let $C$ be a cyclotomic polynomial. Is the set $N(C)$ infinite ?

\end{question}

\section{Proof of Proposition 1 - first part}\label{Prop 1}

In this section and the next one,  we prove Proposition 1 of the introduction.
Let $m, r$ be integers with $m \geqslant 3$ and $r \geqslant 1$, and let $C = \Phi_m^r$.  Assume that ${\rm deg}(C) \leqslant 20$. 

\begin{prop}\label{proposition 1} There exists an automorphism $a : X \to X$ of a projective $K3$ surface $X$ such that the characteristic polynomial of
the restriction of $a^*$ to $T_X$ is equal to $C$.

\end{prop}

The proof of the proposition is divided into several parts, according to the value of $m$. Note first that if $m = p^k$ for some prime
number $p \not = 2$, then Proposition \ref{proposition 1} follows from Proposition \ref{Kondo} :

\begin{prop}\label{p} Suppose that $m = p^k$ where $p$ is a prime number, $p \not = 2$, and $k \geqslant 1$ is an integer. 
Then there exists an automorphism $a : X \to X$ of a  projective $K3$ surface $X$ such that the characteristic polynomial of
the restriction of $a^*$ to $T_X$ is equal to $C$.

\end{prop}

\noindent
{\bf Proof.} We have $C(1) = p^r$ and $C(-1) = 1$, hence Proposition \ref{Kondo} 
implies the existence of an automorphism $a : X \to X$ of a  projective $K3$ surface $X$ with the required
properties.

\begin{prop}\label{2p Kondo} Suppose that $m = 2p^k$ where $p$ is a prime number, and $k \geqslant 1$ is an integer. 
Suppose that $r$ is even, and that ${\rm deg}(C)    \equiv \ 4  \ {\rm (mod \ 8)}$.
Then there exists an automorphism $a : X \to X$ of a  projective $K3$ surface $X$ such that the characteristic polynomial of
the restriction of $a^*$ to $T_X$ is equal to $C$.

\end{prop}

\noindent
{\bf Proof.} We have $C(1) = 1$ and $C(-1) = p^r$; since $r$ is even, $C(-1)$ is a square, hence
the conditions of Proposition \ref{Kondo} are satisfied; therefore this
implies the existence of an automorphism $a : X \to X$ of a  projective $K3$ surface $X$ with the required
properties.

\medskip

In the remaining cases, the proofs use modified versions of the results of the previous sections. 
The following lemma is based on results of McMullen in \cite{Mc3}, and will be used in the proof of Proposition \ref{proposition 1}.  Recall from  \cite{Mc3}, \S 2 that an isometry
of a hyperbolic lattice is said to be {\it positive} if it stabilizes a chamber. 

\begin{lemma}\label{to apply McMullen} Let $(L,q)$ be an even unimodular lattice of signature $(3,19)$, and let $t : L \to L$ be an isometry of $L$. Let
$L_1$ and $L_2$ be mutually orthogonal sublattices of $L$ such that $L_1 \oplus L_2$ is of finite index in $L$, that $t(L_1) = L_1$, $t(L_2) = L_2$, 
that the signature of $L_1$ is $(2,{\rm rank}(L_1) - 2)$ and the signature of $L_2$ is $(1,{\rm rank}(L_2) - 1)$. Suppose
moreover that the restriction of $t$ to $L_2$ preserves a connected component of $\{x \in L_2 \otimes_{\bf Z}{\bf R} \ | \ q(x,x) > 0\}$.
Then we have

\medskip
{\rm (i)} The lattice $L$ has an isometry $t' : L \to L$ such that the restriction of $t'$ to $L_2$ is positive, and that $t'$ and $t$ coincide on $L_1$.

\medskip 
{\rm (ii)} 
Let $V \subset L_1 \otimes_{\bf Z} {\bf R}$ be a $2$-dimensional subspace of signature $(2,0)$ and  stable by $t$ such that the 
intersection of $L$ with the orthogonal of $V$ is equal to $L_2$ and that the restriction of $t$ to $V$ is in ${\rm SO}(V)$. 
Then there exists a projective $K3$ surface $X$ and
an automorphism $a : X \to X$ such that $T_X \simeq L_1$, $S_X \simeq L_2$, and  $a^*|T_X = t$. 

\end{lemma} 

\noindent
{\bf Proof.} (i) Set $t_1 = t|L_1$ and $t_2 = t|L_2$. For $i = 1,2$, set $\overline L_i = L_i^{\sharp}/L_i$, and let $\overline q_i$ and $\overline t_i$ be the induced symmetric bilinear forms and isometries; since $L$ is unimodular, 
we have $(\overline L_1,\overline q_1,\overline t_1) \simeq (\overline L_2,- \overline q_2,\overline t_2)$. 
If $L_2$ has no roots, then $t_2$ is a positive isometry in the sense of McMullen \cite{Mc3}, \S 2; otherwise, let $\rho$ be an element of the Weyl group of $L_2$
such that $\rho \circ t_2$ is positive. Set $t'_2 = t_2$ in the first case, and $t_2' = \rho \circ t_2$ in the second one. The elements of the Weyl group of $L_2$ inducent
the identity on $\overline L_2$, hence we have $(\overline L_1,\overline q_1,\overline t_1) \simeq (\overline L_2,- \overline q_2,\overline t'_2)$. This implies that
there exists an isometry $t' : L \to L$ such that $t' | L_1 = t_1$ and $t' | L_2 = t_2'$; the isometry $t'_2$ is positive, and this implies (i). 

\medskip  (ii) 
Applying \cite{Mc3}, Theorem 6.1 to the isometry $t' : L \to L$ constructed in part (i), we conclude that there exists a projective $K3$ surface $X$ with $T_X \simeq L_1$, $S_X \simeq L_2$, and
an automorphism $a : X \to X$ such that $a^* = t'$. By construction, we have $t'|L_1 = t|L_1$, hence the restriction of $a^*$ to $T_X$ is equal to $t|L_1$. This concludes
the proof of the lemma.

\begin{prop}\label{8} Suppose that $C(1) = C(-1) = 1$.
Then there exists an automorphism $a : X \to X$ of a  projective $K3$ surface $X$ such that the characteristic polynomial of
the restriction of $a^*$ to $T_X$ is equal to $C$.

\end{prop}

\noindent
{\bf Proof.} If ${\rm deg}(C)    \equiv \ 4  \ {\rm (mod \ 8)}$, then this follows from Proposition \ref{Kondo}. Suppose that
${\rm deg}(C)    \equiv \ 0  \ {\rm (mod \ 8)}$, and that $m \not = 30,  60$. Then we have
$m = 15, 20, 24$ and $r = 1$ or $2$, or $m = 40, 48$ and $r = 1$. By Proposition \ref{15}  the lattice $\Lambda_{3,19}$ has
an isometry $t$ such that the characteristic polynomial of $t$ is divisible by $C$ and by $(x-1)^4$, and that the
signature of ${\rm Ker}(C(t))$ is $(2, {\rm deg}(C) - 2)$. The same property holds for $m = 30$ and $r = 1$ by Proposition \ref{30}.
If $m = 60$, then by Example \ref{for 60 cyclotomic}, the lattice $\Lambda_{3,19}$ has
an isometry $t$ such that the characteristic polynomial of $t$ is $C S_2$; if $m = 30$ and $r = 2$, then this holds for $C S_1$ by Example \ref{30}.

\medskip
Set  $L_1 = {\rm Ker}(C(t))$, and
let $L_2$ be the orthogonal complement of $L_1$ in $L$. The 
hypotheses of Lemma \ref{to apply McMullen} are fulfilled; hence 
by Lemma \ref{to apply McMullen} there exists a projective $K3$ surface $X$ with $T_X = L_1$, and an automorphism  $a : X \to X$ such that  the characteristic polynomial of
the restriction of $a^*$ to $T_X$ is equal to $C$.

\begin{prop}\label{2p first part} Let $p$ be a prime number, let  $r,k \geqslant 0$ be integers, and let $C = \Phi_{2p^k}^r$. 
Suppose that  ${\rm deg}(C) \leqslant 16$, or $r$ is even
and ${\rm deg}(C) \leqslant 20$. Then there exists an automorphism $a : X \to X$ of a  projective $K3$ surface $X$ such that the characteristic polynomial of
the restriction of $a^*$ to $T_X$ is equal to $C$.

\end{prop}

\noindent
{\bf Proof.} The hypothesis implies that $C(1) = 1$ (if $p \not = 2$) or $C(1) = 2^r$ (if $p = 2$) and $C(-1) = p^r$.  If $r$ is odd, set $C'(x) = (x+1)^2(x-1)^{20 - {\rm deg}(C)}$. If $C = \Phi_5^4$, set $C = S_0$ (cf. Example \ref{Sa}). Then $CC'$ satisfies condition (C 1) and $G_{CC'} = 0$, hence by
Theorem \ref{signature map theorem} the lattice $\Lambda_{3,19}$ has an isometry with characteristic polynomial $CC'$ such
that the signature of ${\rm Ker}(C(t))$ is $(2, {\rm deg}(C) - 2)$.  
If $r$ is even and $p = 2$ or ${\rm deg}(C)    \equiv \ 4  \ {\rm (mod \ 8)}$, then the existence of such an isometry (with $C'$ a power of $x-1$) follows from
Proposition \ref{Kondo}.  We conclude as in the proof of Proposition \ref{8}.

\section{Proof of Proposition 1 - continued}\label{Prop 1 2p}

The aim of this section is to prove Proposition 1 (that is, Proposition \ref{proposition 1})  in the remaining cases; the results are stated in a more general setting
than needed.

\begin{notation} Let $q$ be a prime number. If $V = (V,b)$ is a quadratic form over ${\bf Q}_q$, we denote by 
$d(V) \in {\bf Q}_q^{\times}/{\bf Q}_q^{\times 2}$ its determinant, and by $w(V) \in {\rm Br}_2({\bf Q}_q)$ its Hasse-Witt
invariant. 

\end{notation}

\begin{lemma}\label{2-adic} Let $V$ be a quadratic form over ${\bf Q}_2$. Then $V$ contains an even, unimodular ${\bf Z}_2$-lattice
if and only if

\medskip
$\bullet$ ${\rm dim}(V)  \equiv 2  \ {\rm (mod \ 4)}$, $d(V) = -1$ and $w(V) = 0$ or $d(V) = 3$ and $d(V) = 1$;

\medskip
$\bullet$ ${\rm dim}(V)  \equiv \ 0  \ {\rm (mod \ 4)}$, $d(V) = 1$ and $w(V) = 1$ or $d(V) = 5$ and $d(V) = 0$.

\end{lemma}

\noindent
{\bf Proof.} Let $H = \langle 1,-1 \rangle$ and $N = \langle 2,6 \rangle$. By \cite{bg}, Proposition 5.2, we see that $V$ 
contains an even, unimodular ${\bf Z}_2$-lattice
if and only if
$V$ is an orthogonal sum of copies of $H$ and $N$; the lemma follows by computing the invariants of these orthogonal sums.

\begin{notation} Let $K$ be a field, and let $E$ be an \'etale $K$-algebra with a $K$-linear involution $\sigma : E \to E$; set
$E_0 = \{x \in E \ | \sigma(x) = x \}$. Let $\lambda \in E_0^{\times}$. We denote by $b_{\lambda}$ the quadratic form
$b_{\lambda} : E \times E \to K$ given by $b_{\lambda}(x,q) = {\rm Tr}_{E/K}(\lambda x \sigma y)$. 

\end{notation}

\begin{prop}\label{18 20} Let $p$ be a prime number, $p \not = 2$, let  $r,k \geqslant 0$ be integers, and let $C = \Phi_{2p^k}^r$. 
Suppose that if $r$ is even, then ${\rm deg}(C)    \equiv \ 4  \ {\rm (mod \ 8)}$. The lattice
$(L,q) = \Lambda_{3,19}$ has an isometry $t$ with characteristic polynomial $C C'$, where $C'(x) = (x^2-1)(x-1)^{20 - {\rm deg}(C)}$,
such that the sublattice ${\rm Ker}(C(t))$ has signature $(2,18)$ and that the restriction of $t$ to ${\rm Ker}(C'(t))$ stabilizes
one of the connected components of 
$\{x \in {\rm Ker}(C'(t)) \otimes_{\bf Z} {\bf R} | \  q(x,x) > 0 \}$.

\end{prop}

\noindent
{\bf Proof.} Set $2n = 22 - {\rm deg}(C)$, and let $U$ be the $\bf Q$-vector space with basis $e_1,\dots,e_n$, $f_1,\dots,f_n$; let
$Q: U \times U \to {\bf Q}$ be the orthogonal sum of the quadratic form equal to $\langle 2, -2p^r \rangle$ on the subspace generated by $e_1$ and $f_1$, and of the diagonal form $\langle -2,\dots, -2 \rangle$ on the subspace generated by $e_i, f_j$ for $i,j \not = 1$. Let $T : U \to U$
be the isometry of $Q$ given by $T(f_1) = -f_1$ and by $T(e_i = e_i$ for all $i$, $T(f_i) = f_i$ if $i \not = 1$. 

\medskip Suppose first that $r$ is even. Since $C(1) = 1$, $C(-1) = p^r$ and ${\rm deg}(C)  \equiv \ 0  \ {\rm (mod \ 4)}$, the polynomial
satisfies condition (C 1). Moreover, $G_F = 0$, since $C$ is a power of an irreducible polynomial. 
We are assuming that ${\rm deg}(C)    \equiv \ 4  \ {\rm (mod \ 8)}$, hence ${\rm deg}(C) -2   \equiv \ 2  \ {\rm (mod \ 8)}$. Therefore Theorem 
\ref{signature map theorem} implies that $\Lambda_{2,{\rm deg}(C) - 2}$ has an isometry $T'$ with characteristic polynomial $C$. Note that
$(U,Q)$ contains a lattice isomorphic to $\Lambda_{1,2n-1}$ stable by $T$, hence $\Lambda_{3,19}$ has an isometry with the required
properties.

\medskip Assume now that $r$ is odd. Set $F = {\bf Q}[x]/(\Phi_{2p^k})$, and let $\alpha \in F$ be the image of $x$. Let $\sigma_F : F \to F$ be
the involution induced by $\alpha \mapsto \alpha^{-1}$, and let $F_0$ be the fixed field of this involution. Let $E_0$ be an
extension of $F_0$ of degree $r$ that is linearly independent of $F$ over $F_0$, and set $E = E_0 \otimes_{F_0} F$. 
Then $E$ is a field, and the characteristic polynomial of $\alpha \in E$ is equal to $C = \Phi_{2p^k}^r$. Let $\sigma$ be the extension
of $\sigma_F$ to $E$. 

\medskip
If $q$ is a prime number, set $E_q = E \otimes _{\bf Q} {\bf Q}_q$ and $(E_q)_0 = E_0 \otimes _{\bf Q} {\bf Q}_q$. With the notation
of \cite{BT},
let $\lambda_p \in (E_p)_0^{\times}$ be such that $\partial (E_p, b_{\lambda_p}, \alpha) = - \partial(U,Q,T)$; if $q \not = p$, let
$\lambda_q \in (E_q)_0^{\times}$ be such that $\partial ((E_q, b_{\lambda_q}, \alpha) = 0$ and that $(E_q, b_{\lambda_q})$
contains an even, unimodular ${\bf Z}_q$-lattice stabilized by $\alpha$; this is possible by \cite{BT}, Proposition 7.1, Proposition 9.1
and the fact that ${\rm det}(E_q, b_{\lambda_q}) = p$ and ${\rm deg}(C) \equiv \ p-1  \ {\rm (mod \ 4)}$. Let $\lambda_{\infty} \in {\bf R}^{\times}$
be such that the signature of $(E \otimes_{\bf Q} {\bf R}, b_{\lambda_{\infty}})$ is equal to $(2,{\rm deg}(C) - 2)$. 

\medskip
For all prime numbers $q$, set $U_q = (U,Q) \otimes _{\bf Q} {\bf Q}_q$ and $W_q = (E_q,b_{\lambda_q})$; we have
$d(U_q) = p$ and $d(W_q) = -p$. Note that this implies that $w(W_q \oplus U_q) = w(W_q) + W(U_q)$. If $q \not = 2,p$, we have
$w(W_q) = w(U_q) = 0$. 

\medskip
Set $W_{\infty} = (E \otimes_{\bf Q} {\bf R}, b_{\lambda_{\infty}})$, and set $w(W_{\infty}) = w_2(W_{\infty})$ in ${\rm Br}_2(\bf R)$. We have
$w(W_{\infty}) = 0 \iff p \equiv \ 3  \ {\rm (mod \ 4)} \iff n \equiv \ 0  \ {\rm (mod \ 2)}$ and $w(W_{\infty}) = 1 \iff p \equiv \ 1  \ {\rm (mod \ 4)}
 \iff n \equiv \ 1  \ {\rm (mod \ 2)}$.
 
 \medskip
 By construction, $W_p \oplus U_p$ contains a unimodular lattice, hence $w(W_p) = w(U_p)$, and we have $w(U_p) = 0 \iff 
 p  \equiv \ \pm 1 \ {\rm (mod \ 8)}$.
 
 \medskip Together with Lemma \ref{2-adic}, this allows us to compute $w(W_q)$ for all $q$, as follows.
 
 \medskip Assume first that $p \equiv \ 3  \ {\rm (mod \ 4)}$. Then $w(W_{\infty}) = 0$ and 
 $$w(W_p) = 0 \iff p  \equiv \ \ 7 \ {\rm (mod \ 8)}.$$
 By Lemma \ref{2-adic}, we
 have 
  $$w(W_2) = 0 \iff p  \equiv \ \ 7 \ {\rm (mod \ 8)}.$$ Since $w(W_q) = 0$ if $q \not = 2,p$, the sum of the invariants $w(W_q)$ (for $q$ 
  a prime number) and $w(W_{\infty})$ is 0.
  
  \medskip Suppose that $p \equiv \ 1  \ {\rm (mod \ 4)}$. Then $w(W_{\infty}) = 1$ and 
 $$w(W_p) = 0 \iff p  \equiv \ \ 1 \ {\rm (mod \ 8)}.$$
 By Lemma \ref{2-adic}, we
 have 
  $$w(W_2) = 0 \iff p  \equiv \ \ 5 \ {\rm (mod \ 8)}.$$ Again, since $w(W_q) = 0$ if $q \not = 2,p$, the sum of the invariants $w(W_q)$ (for $q$ 
  a prime number) and $w(W_{\infty})$ is 0.

\medskip
We have $w_2(E_q,b_{\lambda_q}) = w(b_1) + {\rm cor}_{E_q/{\bf Q}_q}
(\lambda_q,d)$ for all prime numbers $q$ (see \cite {B 22}, Proposition 5.4).  Let $\mathcal V$ be the set of all places of $\bf Q$; since $b_1$ is a global form, the above argument shows that
$\underset {v \in \mathcal V} \sum  {\rm cor}_{E_v/{\bf Q}_v}(\lambda_p,d) = 0$.
By \cite{B 22}, Theorem 9.6, this implies that there exists $\lambda \in E_0^{\times}$ such
that 
$(E,b_{\lambda}) \otimes_{\bf Q} {\bf Q}_q \simeq (E_p,b_{\lambda_p})$ for all $q$.

\medskip
Let $(V,B,t)$ be the orthogonal sum of $(E,b_{\lambda},\alpha)$ and $(U,Q,T)$. 
Set $V_2 = (V,B) \otimes_{\bf Q}{\bf Q}_2$. We have $d(V) = -1$, and  $w(V_2) = w(W_2) + w(U_2)$.
Recall that ${\rm deg}(C) \equiv \ p-1  \ {\rm (mod \ 4)}$, hence ${\rm dim}(W_2) \equiv \ p-1  \ {\rm (mod \ 4)}$; we
have ${\rm dim}(U) = 22 - {\rm deg}(C)$, hence ${\rm dim}(U)  \equiv \ p + 1  \ {\rm (mod \ 4)}$. This implies that

\medskip
If $p \equiv \ 1  \ {\rm (mod \ 4)}$, then 

\medskip

\centerline {$w(U_2) = 0 \iff p \equiv \ 1  \ {\rm (mod \ 8)}$ and $w(W_2) = 0 \iff  p \equiv \ 5  \ {\rm (mod \ 8)}$;}

\medskip
If $p \equiv \ 3  \ {\rm (mod \ 4)}$, then 

\medskip

\centerline {$w(U_2) = 0 \iff p \equiv \ 3  \ {\rm (mod \ 8)}$ and $w(W_2) = 0 \iff  p \equiv \ 7  \ {\rm (mod \ 8)}$.}

\medskip
In both cases, we have $w(W_2) + w(U_2) = 1$, hence $w(V_2) = 1$. 

\medskip The quadratic form $V$ has determinant $-1$, signature $(3,19)$, $w(V_2) = 1$ and all the other Hasse-Witt
invariants of $V$ are trivial. This implies that $V$ is isomorphic to $\Lambda_{3,19} \otimes_{\bf Z} {\bf Q}$. The
characteristic polynomial of $t$ is $CC'$. The quadratic form $V$ contains an
even unimodular lattice stabilized by $t$ everywhere locally; this is clear by construction at all prime numbers $q \not = 2$, and
for $q = 2$ it follows from Takada, \cite{T}, Theorem 4.2.
The intersection of these lattices
is an even unimodular lattice stabilized by the isometry $t$; this lattice is isomorphic to $\Lambda_{3,19}$.  By construction,
$t$ has the required properties. This concludes the proof of the Proposition.

\begin{prop}\label{2 cyclotomic} Let $C = \Phi_{2^k}^r$ with $k = 2$ and $r = 9$ or $k = 3$ and $r = 5$. . Then the lattice
$(L,q) = \Lambda_{3,19}$ has an isometry $t$ with characteristic polynomial $C C'$, where $C'(x) = (x^2-1)(x-1)^{20 - {\rm deg}(C)}$,
such that the sublattice ${\rm Ker}(C(t))$ has signature $(2,18)$ and that the restriction of $t$ to ${\rm Ker}(C'(t))$ stabilizes
one of the connected components of 
$\{x \in {\rm Ker}(C'(t)) \otimes_{\bf Z} {\bf R} | \  q(x,x) > 0 \}$.

\end{prop}

\noindent
{\bf Proof.} Set $E = {\bf Q}/[x]/(\Phi_{2^k})$, with $k = 2$ or $k = 4$. We denote by $x \mapsto \overline x$ the complex conjugation,
and let $E_0$ be the fixed  subfield of $E$ : we have $E_0 = {\bf Q}$ if $r = 2$ and 
$E_0 = {\bf Q}(\sqrt 2)$ if $k = 4$. 

\medskip
Suppose first that $k = 2$, and let $X = (X,q_X,t_X)$ 
be defined by $X = E$, $q_X(x,y) = {\rm Tr}_{E/{\bf Q}}({1 \over 2}x \overline y)$; let $t_X$ be induced by multiplication by $i = \zeta_4$; note
that $t_X$ is an isometry of $q_X$ with characteristic polynomial $\Phi_4$.
Let $W_2 = (W_2,q_2,t_2)$  be the orthogonal sum of a copy of $X$ with 9 copies of $-X$. 
We have ${\rm dim}(W_2) = 18$, $d(W_2) = 1$, and $w(W_2) = 0$.  The signature of $W_2$ is $(2,16)$, and the characteristic polynomial of $t_2$ is $\Phi_4^9$. 

\medskip
Let $U_2$ be the ${\bf Q}$-vector space of basis $e_1,e_2,f_1,f_2$ and $q_2 : U_2 \times U_2 \to {\bf Q}$ be the
quadratic form such that $q_2(e_1,e_1) = 1$ and $q_2(f_1,f_1) = q_2(e_2,e_2) = q_2(f_2,f_2) = -1$. 
Let $t_2 : U_2 \to U_2$ be the isometry given by $t_2(f_2) = -f_2$ and $t_2(e_i) = e_i$ for $i = 1,2$, $t_2(f_2) = f_2$. 
We have  ${\rm dim}(U_2) = 4$, $d(U_2) =  -1$,  $w(U_2) = 1$ at $2$ and $\infty$, and $0$ elsewhere. 

\medskip
If $k = 2$, we set $(V,q,t) = (W_2,q_2,t_2) \oplus (U_2,q_2,t_2)$. The signature of $V$ is $(3,19)$, and $d(V) = -1$, $w(V) = 1$
at $2$ and $\infty$, and $0$ elsewhere. The characteristic polynomial of $t$ is $\Phi_4^9(x)(x^2-1)(x-1)^2$. 

\medskip
Assume now that $k = 4$. Let  $X = (X,q_X,t_X)$ 
be defined by $X = E$, $q_X(x,y) = {\rm Tr}_{E/{\bf Q}}({\sqrt 2 \over 4}x \overline y)$; let $t_X$ be induced by multiplication by $ \zeta_8$; note
that $t_X$ is an isometry of $q_X$ with characteristic polynomial $\Phi_8$. Let $Y = Y,q_Y,t_Y)$ be defined by 
$Y = E$, $q_Y(x,y) = {\rm Tr}_{E/{\bf Q}}({1 \over 4}x \overline y)$; let $t_Y$ be induced by multiplication by $ \zeta_8$. Let
$W_4 = (W_4,q_4,t_4)$ be the orthogonal sum of $X$ with $4$ copies of $-Y$. We have ${\rm dim}(W_4) = 20$, $d(W_2) = 1$, and $w(W_2) = 1$
at $2$ and at  $\infty$, and $0$ elsewhere. 
The signature of $W_2$ is $(2,18)$, and the characteristic polynomial of $t_2$ is $\Phi_8^5$. 

\medskip
Let $U_2$ be the ${\bf Q}$-vector space of basis $e_1,f_1$ and $q_2 : U_2 \times U_2 \to {\bf Q}$ be the
quadratic form such that $q_2(e_1,e_1) = 1$ and $q_2(f_1,f_1) = -1$. 
Let $t_2 : U_2 \to U_2$ be the isometry given by $t_2(f_2) = -f_2$ and $t_2(e_1) = e_1$. 
We have  ${\rm dim}(U_2) = 2$, $d(U_2) =  -1$,  $w(U_2) = 0$.

\medskip
If $k = 2$, we set $(V,q,t) = (W_2,q_2,t_2) \oplus (U_2,q_2,t_2)$. The signature of $V$ is $(3,19)$, and $d(V) = -1$, $w(V) = 1$
at $2$ and $\infty$, and $0$ elsewhere. The characteristic polynomial of $t$ is $\Phi_4^9(x)(x^2-1(x-1)^2)$. 

\medskip
If $k = 4$, we set $(V,q,t) = (W_4,q_4,t_4) \oplus (U_4,q_4,t_4)$. The signature of $V$ is $(3,19)$, and $d(V) = -1$, $w(V) = 1$
at $2$ and $\infty$, and $0$ elsewhere. The characteristic polynomial of $t$ is $\Phi_8^5(x)(x^2-1)$. 

\medskip
In both cases, $(V,q,t)$ contains an even unimodular lattice stabilized by $t$ everywhere locally; at the prime $2$, this
follows from Takada 
\cite{T}, Theorem 4.2.  Let $L$ be the intersection of these lattices; $L$ is stabilized by $t$, and we have $L \simeq \Lambda_{3,19}$. This concludes
the proof of the proposition.

\medskip

\noindent
{\bf Proof of Proposition \ref{proposition 1}.} If $m = p^k$ where $p$ is a prime number with $p \not = 2$, then the proposition
follows from Proposition \ref{p}; if $C(1) = C(-1)$, then from Proposition \ref{8}. Suppose that $C = \Phi_{2p^k}^r$. If $r$ is
even of if ${\rm deg}(C) \leqslant 16$, it is a consequence of Proposition \ref{2p first part}. Suppose that $r$ is odd and
that ${\rm deg}(C) = 18$ or $20$; apply Proposition \ref{18 20} f $p \not = 2$, and Proposition \ref{2 cyclotomic} if $p = 2$.

\bigskip
Proposition \ref{proposition 1} implies the following result :

\begin{coro}\label{OM} Let $m$ be an integer such that $m \geqslant 1$ and that $\varphi(m) \leqslant 20$. Then  there exists an automorphism $a : X \to X$  of a projective $K3$ surface $X$ 
inducing multiplication by a primitive $m$-th root of unity on $T_X$.

\end{coro}

\noindent
{\bf Proof.} For $m = 1$, we can take the identity, and there are many examples of automorphisms of projective $K3$ surfaces $X$ inducing $-{\rm id}$ on $T_X$ (see for instance
\cite{H}, Corollary 15.2.12). Suppose that $m \geqslant 3$, and let $C = \Phi_m$; Proposition \ref{proposition 1} implies that there exists an automorphism $a : X \to X$
of a projective $K3$ surface such that the characteristic polynomial of $a^*|T_X$ is equal to $C$; hence $a^*|T_X$ acts by multiplication
by a primitive $m$-th  root of unity. 

\begin{remark}\label{finite order MO}
Corollary \ref{OM} follows from results of  Machida-Oguiso \cite{MO}, Xiao \cite{X} and Zhang \cite {Z} when $m \not = 60$; more precisely, they prove the existence of an automorphism $a : X \to X$ of {\it finite order}
inducing multiplication by a primitive $m$-th  root of unity on $T_X$. They
 also show that this is not the case for $m = 60$; in
the next section we give another proof of this result.

\end{remark} 

\section{The primitive $60$-th roots of unity}\label{60 section}

The following result was proved by Machida-Oguiso \cite{MO}, Xiao \cite{X} and Zhang \cite {Z}  : 

\begin{prop}\label{no 60}
There does not exist 
any automorphism of  finite order of a projective $K3$ surface inducing multiplication by a primitive $60$-th root of unity on its transcendental lattice. 

\end{prop}

The aim of this section is to give another proof of Proposition \ref{no 60}. Set $C = \Phi_{60}$. 

\begin{prop}\label{60.1} Let $L$ be an even unimodular lattice of signature $(3,19)$. The lattice $L$ does not have any isometry $t : L \to L$ having the
following properties :

\medskip
{\rm (i)} The characteristic polynomial of $t$ is $CC'$, where $C'$ is a product of cyclotomic polynomials.

{\rm (ii)} The signature of the sublattice $L_C = {\rm Ker}(C(t))$ of $L$ is $(2,14)$. 

\end{prop} 

The main ingredient of the proof of Proposition \ref{60.1} is the following lemma~:

\begin{lemma}\label{60.2} Set $C' = \Phi_{12}$, and let $M$ be an even unimodular lattice of signature $(2,18)$. The lattice $M$ does not have any isometry $t : M \to M$  having the
following properties :

\medskip
{\rm (i)} The characteristic polynomial of $t$ is $CC'$.

{\rm (ii)} The signature of the sublattice $M_C = {\rm Ker}(C(t))$ of $M$ is $(2,14)$. 

\end{lemma} 

We give two proofs of this lemma; the first one is based on some results of \cite{B 22}, the second one is a direct proof. 

\medskip 
\noindent
{\bf First proof of Lemma \ref{60.2}.} Set $F = CC'$, and note that $F$ satisfies condition (C 1). By Example \ref{12 and 60}, we have $G_F = {\bf Z}/2{\bf Z}$.  

\medskip Set $I = \{ C,C' \}$. Since  $G_F = {\bf Z}/2{\bf Z}$, we have $C_0(I) = C(I)$. Let
$c : I \to {\bf Z}/2{\bf Z}$ be such that $c(C) = 1$ and $c(C') = 0$, and let $c' : I \to {\bf Z}/2{\bf Z}$ be such that $c'(C) = 0$ and $c'(C') = 1$.

\medskip As in \cite{B 22}, \S 9 and \S 12, we define a homomorphism $\epsilon_F^{\rm finite} : C(I)  \to {\bf Z}/2{\bf Z}$. 

\medskip Let $\tau$ be a signature map with characteristic polynomial $F$ and maximum $(2,18)$ such that $\tau(C) = (2,14)$ and $\tau(C') = (0,4)$, and let 
$\epsilon_{\tau}^{\infty} : C(I) \to {\bf Z}/2{\bf Z}$ be the associated homomorphism (see \cite {B 22}, \S 9 and \S 12). We obtain a homomorphism
$\epsilon_{\tau} : G_F \to {\bf Z}/2{\bf Z}$ by setting $\epsilon_{\tau} = \epsilon_F^{\rm finite} + \epsilon_{\tau}^{\infty}$.
By \cite{B 22}, Theorem 12.1, there
exists an even unimodular lattice $M$ having an isometry $t : M \to M$ with properties (i) and (ii) if and only if $\epsilon_{\tau} = 0$.

\medskip With the notation of \cite{B 22}, \S 9, we have $a^{\infty}_{\tau}(C) = 1$ and $a^{\infty}_{\tau}(C') = 0$. This implies that $\epsilon_{\tau}^{\infty} (c) = 1$ and
$\epsilon_{\tau}^{\infty} (c') = 0$.

\medskip Similarly, let $\tau'$ be a signature map with characteristic polynomial $F$ and maximum $(2,18)$ such that $\tau(C) = (0,16)$ and $\tau(C') = (2,2)$, and let 
$\epsilon_{\tau}^{\infty} : C(I) \to {\bf Z}/2{\bf Z}$ be the associated homomorphism. We have $a^{\infty}_{\tau'}(C) = 0$ and $a^{\infty}_{\tau'}(C') = 1$, hence  $\epsilon_{\tau'}^{\infty} (c) = 0$ and
$\epsilon_{\tau'}^{\infty} (c') = 1$.
We obtain a homomorphism
$\epsilon_{\tau'} : G_F \to {\bf Z}/2{\bf Z}$ by setting $\epsilon_{\tau'} = \epsilon_F^{\rm finite} + \epsilon_{\tau'}^{\infty}$.

\medskip Both $C$ and $C'$ satisfy condition (C 1), and since they are irreducible, we have $G_C = G_{C'} = 0$. Therefore by Theorem \ref{basic theorem} there exists an even unimodular lattice
$N$ of signature $(0,16)$ having an isometry with characteristic polynomial $C$, and an even unimodular lattice $N'$ of signature $(2,2)$ having an isometry
with characteristic polynomial $C'$. The lattice $N \oplus N'$ is even unimodular of signature $(2,18)$, and has an isometry of characteristic polynomial $F$ and
of signature map $\tau'$. Applying \cite{B 22}, Theorem 12.1, this implies that $\epsilon_{\tau'} = 0$.
Since $\epsilon_{\tau'} = \epsilon_F^{\rm finite} + \epsilon_{\tau'}^{\infty}$, we have  $\epsilon^{\rm finite}_F(c) = 0$ and $\epsilon^{\rm finite}_F(c') = 1$.

\medskip On the other hand, we have  $\epsilon_{\tau} = \epsilon_F^{\rm finite} + \epsilon_{\tau}^{\infty}$, and this 
implies that $\epsilon_{\tau} \not = 0$; therefore there does not exist any even unimodular lattice $M$ having an isometry $t : M \to M$ with properties (i) and (ii).

\medskip
The second proof of Lemma \ref{60.2} uses the notion of Hasse-Witt invariant of a quadratic form.

\begin{notation} Let $K$ be a field of characteristic $\not = 2$, let $V$ be a finite dimensional $K$-vector space, and let $q : V \times V \to K$ be a non-degenerate quadratic
form. The Hasse-Witt invariant of $q$ is denoted by $w_2(q)$; it is an element of ${\rm Br}_2(K)$.

\end{notation}

If $K$ is a $p$-adic field or ${\bf R}$, then ${\rm Br}_2(K)$ is a group of order two, that we identify with $\{0,1\}$; see \cite{S}, Chapitre IV  for the properties
of Hasse-Witt invariants of quadratic forms that are needed here.

\medskip 
\noindent
{\bf Second proof of Lemma \ref{60.2}.} Set $F = C C'$. Suppose that the even unimodular lattice $(M,q)$ of signature $(2,18)$ has an isometry $t : M \to M$ with
characteristic polynomial $F$, and let us show that (ii) does not hold.

\medskip Set $M_1 = {\rm Ker}(C(t))$ and let $q_1$ be the restriction of $q$ to $M_1 \times M_1$; similarly, set $M_2 = {\rm Ker}(C'(t))$ and let $q_2$ be the restriction of $q$ to $M_2 \times M_2$.  Set $V = M \otimes_{\bf Z} {\bf Q}$, $V_1 = M_1 \otimes_{\bf Z} {\bf Q}$ and $V_2 = M_2 \otimes_{\bf Z} {\bf Q}$.
Since $C$ and $C'$ are distinct irreducible polynomials, the quadratic space $(V,q)$ is the orthogonal sum of $(V_1,q_1)$ and $(V_2,q_2)$. 

\medskip
The resultant of $C$ and $C'$ is $5^4$. This implies that if $p$ is a prime number with $p \not = 5$, then $(M,q) \otimes_{\bf Z} {\bf Z_p}$ is the orthogonal sum
of $(M_1,q_1) \otimes_{\bf Z} {\bf Z}_p$ and $(M_2,q_2) \otimes_{\bf Z} {\bf Z}_p$; therefore the ${\bf Z_p}$-lattices 
$(M_1,q_1) \otimes_{\bf Z} {\bf Z}_p$ and $(M_2,q_2) \otimes_{\bf Z} {\bf Z}_p$  are unimodular.

\medskip Let $p$ be a prime number $p \not = 2, 5$. Since $(V_1,q_1) \otimes_{\bf Q} {\bf Q}_p$ and 
$(V_2,q_2) \otimes_{\bf Q} {\bf Q}_p$ contain unimodular ${\bf Z}_p$-lattices, we have $w_2(q_1) = w_2(q_2) = 0$ in ${\rm Br}_2({\bf Q}_p)$. 

\medskip Let $U$ be the hyperbolic plane; if $n$ is an integer with $n \geqslant 1$, we denote by $U^n$ the orthogonal sum of $n$ copies of $U$. Since $(M_1,q_1) \otimes_{\bf Z} {\bf Z}_2$ and $(M_2,q_2) \otimes_{\bf Z} {\bf Z}_2$ are unimodular, we have
$(M_1,q_1) \otimes_{\bf Z} {\bf Z}_2 \simeq U^8 \otimes_{\bf Z} {\bf Z}_2$ and $(M_2,q_2) \otimes_{\bf Z} {\bf Z}_2 \simeq  U^2 \otimes_{\bf Z} {\bf Z}_2$ (see for instance \cite{bg}, Proposition 5.2). This
implies that $w_2(q_1) = 0$ and $w_2(q_2) = 1$  in ${\rm Br}_2({\bf Q}_2)$. 

\medskip Let $K$ be the cyclotomic field of the $12$-th roots of unity, and let $K_0$ be its maximal real subfield The prime $5$ is inert in the extension $K_0/{\bf Q}$, and splits in 
$K/K_0$. This implies that there exists a degree $2$ polynomial $f \in {\bf Z}_5[x]$ such that $C' = f f^*$ in $ {\bf Z}_5[x]$, where $f^*(x) = x^2 f(x^{-1})$,
and such that $f \not = f^*$. Let us denote by $t_2$ the restriction of $t$ to $V_2$, and note that ${\rm Ker}(f(t_2))$ is an isotropic subspace of dimension $2$ of $V_2 \otimes_{\bf Q} {\bf Q}_5$; this implies that $(V_2,q_2) \otimes_{\bf Q} {\bf Q}_5 \simeq U^2  \otimes_{\bf Q} {\bf Q}_5$, and therefore $ w_2(q_2) = 0$ in ${\rm Br}_2({\bf Q}_5)$.

\medskip
Suppose that (ii) holds. Then the signature of $(V_1,q_1)$ is $(2,14)$, and the signature of  $(V_2,q_2)$ is $(0,4)$; hence 
$w_2(q_2) = 0$ in ${\rm Br}_2({\bf R})$. 
This leads to a contradiction, since 
$w_2(q_2) = 1$  in ${\rm Br}_2({\bf Q}_2)$, 
and $ w_2(q_2) = 0$ in ${\rm Br}_2({\bf Q}_p)$ for all prime numbers $p$ with $p \not = 2$. 

\medskip
\noindent 
{\bf Proof of Proposition \ref{60.1}.} Let $C'$ be a product of cyclotomic polynomials and let $t : L \to L$ be an isometry with characteristic polynomial $C C'$, and set $F = CC'$. Suppose
first that $C'$ is not divisible by $\Phi_{12}$. Then $C$ and $C'$ are relatively prime over ${\bf Z}$. Indeed, ${\rm deg}(C') = 6$, and if $\Phi$
is a cyclotomic polynomial of degree $\leqslant 6$ such that $\Phi$ is not relatively prime to $C$ over ${\bf Z}$, then $\Phi = \Phi_{12}$; this
follows from the values of the resultants of cyclotomic polynomials, see for instance \cite{Ap}.
Since $C$ and $C'$ are relatively prime over ${\bf Z}$, the lattice
$L$ is the orthogonal sum of the even unimodular
lattices $L_C$ and ${\rm Ker}(C'(t))$. If (ii) holds, then the signature of $L_C$ is $(2,14)$; this contradicts the fact that $L_C$ is an even unimodular lattice. 

\medskip Assume now that $\Phi_{12}$ divides $C'$. The polynomial $F$ satisfies condition (C~1); this implies that $C' = \Phi_{12} C''$ such that the irreducible
factors of $C''$ are in $\{x-1,x+1\}$. Therefore $C''$ is relatively prime over ${\bf Z}$ to $C \Phi_{12}$. Set $M = {\rm Ker}(C \Phi_{12}(t))$ and $M' = {\rm Ker}(C''(t))$.
The lattice $L$ is the orthogonal sum of $M$ and $M'$, hence both these lattices are even and unimodular. This implies that the signature of $M'$ is $(1,1)$, and
hence the signature of $M$ is $(2,14)$. By Lemma \ref{60.2}, this is impossible. This completes the proof of the proposition.

\medskip
\noindent 
{\bf Proof of Proposition \ref{no 60}.} Let $a : X \to X$ be an automorphism of a projective $K3$ surface such that $a^*$ induces multiplication by a primitive $60$-th
root of unity on $T_X$. This implies that the characteristic polynomial of $a^*|T_X$ is equal to $\Phi_{60}$. Assume that $a$ is of finite order. Then the characteristic
polynomial of $a^*|S_X$ is a product of cyclotomic polynomials. Since $T_X = {\rm Ker}(\Phi_{60}(a^*))$, the signature of the lattice $ {\rm Ker}(\Phi_{60}(a^*))$ is
$(2,14)$. Proposition \ref{60.1} implies that this is impossible, hence no such automorphism exists.

\section{Two cyclic groups}\label{cyclic}

Recall from the introduction that if $X$ is a projective $K3$ surface, we have the exact sequences
$$1 \to {\rm Aut}_s(X) \to {\rm Aut}(X) \to M_X \to 1,$$ and
$$1 \to N_X \to {\rm Aut}(X) \to {\rm O}(S_X),$$

\medskip
\noindent
where $M_X$ and $N_X$ are finite cyclic groups, of order $m_X$, respectively $n_X$. The group $N_X$ can be identified to a subgroup of $M_X$, hence $n_X$ divides $m_X$.

\medskip The question of determining the possible values of $m_X$ and $n_X$ was raised by Huybrechts in \cite{H}, page 336; the characterization of the pairs
$(m_X,n_X)$ is also of interest. 

\medskip The values of $n_X$ were studied much earlier, by Vorontsov  \cite{V}, Kondo \cite{Ko} and Oguiso-Zhang \cite{OZ}; assuming that ${\rm rank}(T_X) = \varphi(n_X)$,
they give a complete list of these values (see Corollary \ref{unimodular} below). 

\medskip
We start with the integers $m_X$. The following is Corollary 1 of the introduction : 

\begin{coro}\label{corollary 1} Let $m \geq 2$ be an even integer such that $\varphi(m) \leqslant 20$. Then  there exists a projective $K3$ surface $X$ 
with $m_X = m$.

\end{coro}

\noindent
{\bf Proof.} There exist projective $K3$ surfaces $X$ with $m_X=2$, see for instance \cite{H}, Corollary 15.2.12. Suppose that $m \geqslant 4$. By
Proposition \ref{proposition 1} there exists a projective $K3$ surface $X$ and an automorphism $a : X \to X$ such that the characteristic polynomial of the restriction of $a^*$ to $T_X$ is $\Phi_m$. 
Since $m$ is even, we have $m_X = m$.

\begin{remark}\label{OM list} In \cite{MO}, Machida and Oguiso consider a related problem; they are interested in the images of the {\it finite subgroups}
of ${\rm Aut}(X)$ in $M_X$. Their results imply that there exist projective $K3$ surfaces $X$
with $m_X = 28, 30, 32, 34, 38, 40$, $42, 44, 48, 54$ and $66$ (see \cite{MO}, Proposition 4).

\end{remark}





\medskip The possible values of $n_X$ can be deduced from Proposition \ref{Kondo}. As we will see, it is enough to consider
$K3$ surfaces $X$ such that ${\rm rank}(T_X) = \varphi(n_X)$ or ${\rm rank}(T_X) = 2 \varphi(n_X)$.

\begin{prop}\label{1 or 2}  Let $m$ be an integer with $m \geqslant 1$ and  $\varphi(m) \leqslant 20$. The following are equivalent

\medskip
{\rm (i)} There exists a projective $K3$ surface $X$ 
such that $n_X = m$.

\medskip
{\rm (ii)} There exists a projective $K3$ surface $X$ 
such that $n_X = m$ and ${\rm rank}(T_X) = \varphi(m)$ or ${\rm rank}(T_X) = 2 \varphi(m)$.

\end{prop}

This will be proved at the end of the section. 

\medskip
We start with the case where $C$ is a cyclotomic polynomial, i.e. when the rank of $T_X$
is $\varphi(n_X)$.

\begin{coro}\label{VK}  Let $m$ be an integer such that $m \geqslant 3$ and that $\varphi(m) \leqslant 20$; set $C = \Phi_m$. There exists a projective $K3$ surface $X$ such that $n_X = m$ and the rank of $T_X$ is $\varphi(n_X)$ if and only if the following conditions hold :

\medskip
{\rm (i)} $C(-1) = 1$.

\medskip
{\rm (ii)} If $C(1) = 1$, then  $m  \equiv \ 0 \ {\rm (mod \ 2)}$ and ${\rm deg}(C)  \equiv \ 4  \ {\rm (mod \ 8)}$. 

\medskip  The $K3$ surface is unique up to isomorphism. Moreover, the lattice $T_X$ is unimodular if and only if $C(1) = 1$, and
$p$-elementary with $p = C(1)$ if $C(1) > 1$. 

\end{coro}

\noindent
{\bf Proof.} Suppose that conditions (i) and (ii) hold. Then by Proposition \ref{Kondo}, there exists an automorphism $a : X \to X$ of a projective $K3$ surface such
that the restriction of $a^*$ is the identity and that the characteristic polynomial of $a^*|T_X$ is equal to $C$; this implies that  $m$ divides $n_X$.  Since $C = \Phi_m$,
the rank of $T_X$ is equal to $\varphi(m)$; hence $\varphi(m) = \varphi(n_X)$. Suppose that $C(1) > 1$; then $m = p^r$, where $p$ is an odd prime number (cf. Lemma \ref{1 and -1}). Then
either $n_X = m$ or $n_X = 2m$; 
but $\Phi_{2p^r}$
does not satisfy (i), hence $n_X \not = 2p^r$, and this implies that $n_X = m$. Assume now that $C(1) = 1$. By (ii), we have $m  \equiv \ 0 \ {\rm (mod \ 2)}$; since
$m$ divides $n_X$ and $\varphi(m) = \varphi(n_X)$, this implies that $n_X = m$. 

\medskip Conversely, suppose that there exists a projective $K3$ surface $X$ such that $n_X = m$ and ${\rm rank}(T_X) = \varphi(m)$, and let $a$ be a generator of the cyclic group $N_X$. Then 
the restriction of $a^*$ to $S_X$ is the identity and the characteristic polynomial of $a^*|T_X$ is equal to $C$. Therefore Proposition \ref{Kondo} (i) implies that $C(-1)$ is a square;
since $C$ is a cyclotomic polynomial, we have
$C(-1) = 1$, hence  (i) holds. Proposition \ref{Kondo} (ii) implies that if  $C(1) = 1$, then  ${\rm deg}(C)  \equiv \ 4  \ {\rm (mod \ 8)}$. Moreover, the hypothesis $C(1) = 1$ implies that $T_X$ is unimodular; this implies that $n_X  \equiv \ 0 \ {\rm (mod \ 2)}$, hence $m  \equiv \ 0 \ {\rm (mod \ 2)}$, and therefore (ii) holds.

\medskip The uniqueness of the $K3$-surface follows from Proposition \ref{Kondo}. 
 If $C(1) = 1$ then $T_X$ is unimodular (cf. Proposition \ref{Kondo}). Conversely, suppose that $T_X$ is unimodular; then $C$ satisfies
condition (C~1). Since $C$ is a cyclotomic polynomial, this implies that $C(1) = 1$. If $C(1) > 1$, then by 
Proposition \ref{p-primary} the lattice $T_X$ is $p$-elementary with $p = C(1)$. Conversely, if $T_X$ is $p$-elementary, then by Proposition
\ref{Kondo} we have $C(1) > 1$.

\bigskip
Using Corollary \ref{VK}, we recover the lists of Vorontsov \cite{V} and Kondo \cite{Ko}, as follows. Set

\medskip
 $A = \{12, \ 28, \ 36, \ 42, \ 44, \ 66 \}$ and $B = \{ 3, \ 9, \ 27, \ 5, \ 25, \ 7, \ 11, \ 13, \ 17, \ 19 \}$. 

\begin{lemma}\label{A and B}
Let $m$ be an integer such that $m \geqslant 3$ and that $\varphi(m) \leqslant 20$. Then we have 

\medskip $m \in A \iff \Phi_m(1) = \Phi_m(-1) = 1$, $m$ is even, and ${\rm deg}(C)  \equiv \ 4  \ {\rm (mod \ 8)}$. 

\medskip $m \in B \iff \Phi_m(-1) = 1$ and $\Phi_m(1) > 1$.

\end{lemma}

\noindent
{\bf Proof.} This follows from Lemma \ref{1 and -1}.

\bigskip Combining Corollary \ref{VK}  and Lemma \ref{A and B}, we obtain 

\begin{coro}\label{unimodular} {\rm (Kondo, Vorontsov)} Let $m \in A \cup B$. Then there exists a projective $K3$ surface $X$ such that
$n_X = m$ and that  the rank of $T_X({X})$ is equal to $\varphi(n_{X})$. Moreover, $X$ is unique up to isomorphism, and

\medskip

\centerline {$n_{X} \in A \iff$
the lattice $T_X$ is unimodular.} 

\medskip

\centerline  {$n_{X} \in B \iff$
the lattice $T_X$ is $p$-elementary with $p = \Phi_{n_X}(1)$.}

 \end{coro}
 
\noindent
{\bf Proof.} See \cite{Ko}, \cite{V}; this also follows from Corollary \ref{VK} and Lemma \ref{A and B}.

\bigskip We now consider the case where ${\rm rank}(T_X) =  2 \varphi(n_X)$.

 \begin{coro}\label{2}  Let $m$ be an integer such that $m \geqslant 3$, that $m  \equiv \ 0 \ {\rm (mod \ 2)}$ and that $\varphi(m) \leqslant 10$. Suppose
 that if $\Phi_m(1) = 1$, then $\varphi(m) \equiv \ 2 \ {\rm (mod \ 4)}$. Then there  exist infinitely many projective $K3$ surfaces $X$ such that $n_X = m$, 
 and that ${\rm rank}(T_X) =  2 \varphi(m)$. 
 

 

 \end{coro}
 
 \noindent
 {\bf Proof.} Set $C = \Phi_m^2$, and note that
 $C$ satisfies conditions (i) and (ii) of Proposition \ref{Kondo}. Therefore there exist infinitely many projective $K3$ surfaces $X$ having automorphisms $a : X \to X$ 
 such that $a^*|S_X$ is the identity, and that the characteristic polynomial of $a^*|T_X$ is $C$. If $m$ is a power of $2$, then $m = n_X$. Otherwise, 
 we have $m = 2m'$ where $m' = 3,7,9$ or $11$; therefore $n_X = m$ or $n_X = 2m = 4m'$. 
 
 \medskip If $Y$ is a $K3$ surface with $n_Y = 4m'$ and ${\rm rank}(T_Y) = 2 \varphi(m) = \varphi(4m')$, then $Y$
 has an automorphism $b : Y \to Y$ such that $b^*|T_Y$ is the identity and the characteristic polynomial of $b^*|T_Y$ is 
 $C' = \Phi_{4m'}$; by Corollary \ref{VK}, the surface $Y$ is unique up to isomorphism. 
 
 \medskip Therefore there exist infinitely many non-isomorphic projective $K3$ surfaces 
 $X$ such that  $n_X = m$. This concludes the proof of the corollary. 
 

\bigskip 
Set $C = \{4,6,8,14,16,18,22\}$.

\begin{lemma}\label{C}  Let $m$ be an integer such that $m \geqslant 3$, and that $\varphi(m) \leqslant 10$. We have

\medskip

$m \in C \iff$ $m$ is even, and if $\Phi_m(1) = 1$, then $\varphi(m) \equiv \ 2 \ {\rm (mod \ 4)}$.

 \end{lemma}
 
 \noindent
 {\bf Proof.} This follows from Lemma \ref{1 and -1}. 
 
 \medskip 
 Corollary \ref{2} and Lemma \ref{C} imply the following
 
 \begin{coro}\label{list C} If $m \in C$, then there  exist infinitely many projective $K3$ surfaces $X$ such that $n_X = m$, 
 and that  ${\rm rank}(T_X) =  2 \varphi(n_X)$.

 \end{coro}

Finally, we show that Corollaries \ref{VK} and \ref{2} characterize the possible values of $n_X$. We start with a lemma :

 \begin{lemma}\label{r} Let $m, r$ be integers with $m \geqslant 3$ and $r \geqslant 1$, and let $C = \Phi_m^r$. Conditions {\rm (i)} and {\rm (ii)}
 of Proposition \ref{Kondo} hold
 for $C$ if and only if 
 
 \medskip
 {\rm (i)} $r  \equiv \ 1 \ {\rm (mod \ 2)}$, $\Phi_m(-1) = 1$ and if $\Phi_m(1) = 1$, then $\varphi(m) \equiv \ 4 \ {\rm (mod \ 8)}$.
 
 \medskip
 {\rm (ii)} $r  \equiv \ 2 \ {\rm (mod \ 4)}$, and if $\Phi_m(1) = 1$, then $\varphi(m) \equiv \ 2 \ {\rm (mod \ 4)}$.
  
 \medskip
 {\rm (iii)} $r  \equiv \ 0 \ {\rm (mod \ 4)}$, and $\Phi(1) > 1$.

 \end{lemma}
 
 \noindent
 {\bf Proof.} This is straightforward verification, using Lemma \ref{1 and -1}.

\begin{coro}\label{characterization} Let $m$ be an integer such that $m \geqslant 3$ and that $\varphi(m) \leqslant 20$. There exists a projective $K3$ surface $X$ with $n_X = m$ and ${\rm rank}(T_X) =  r \varphi(n_X)$ if and only if 

\medskip
{\rm (i)} $r = 1$ and $m \in A \cup B$.

\medskip
{\rm (ii)} $r = 2$ and $m \in C \cup \{3,5,7,11\}$.

\medskip
{\rm (iii)} $r = 3$ and $m \in \{3,5,7,9,12\}$.

\medskip
{\rm (iv)} $r = 4$ and $m \in \{3,4,5,8\}$.

\medskip
{\rm (v)} $r = 5$ and $m \in \{3, 5\}$.

\medskip
{\rm (vi)} $r = 6$ or $10$ and $m \in \{3,4,6\}$.

\medskip
{\rm (vii)} $r = 8$ and $m \in \{3,4\}$.

\medskip
{\rm (viii)} $r = 7$ or $9$ and $m = 3$.

\end{coro}

\noindent
{\bf Proof.} 
This follows from Proposition \ref{Kondo} and Lemma \ref{r}.

\begin{coro}\label{summary} Let $m$ be an integer such that $m \geqslant 3$ and that $\varphi(m) \leqslant 20$. There exists a projective $K3$ surface $X$ with $n_X = m$ 
$\iff m \in A \cup B \cup C$. 

\end{coro}

\noindent
{\bf Proof.} This is an immediate consequence of  Corollary  \ref{characterization}. 

\medskip
\noindent
{\bf Proof of Proposition \ref{1 or 2}.} It is clear that (ii) $\implies$ (i). Suppose (i). Then by Corollary \ref{summary}, we have 
$n_X \in A \cup B \cup C$, and (i) follows from Corollary \ref {unimodular} and Corollary \ref {list C}.

\begin{notation}\label{t} If $X$ is a $K3$ surface, set $r_X = {{\rm rank}(T_X) \over \varphi(n_X)}$ and $t_X = {{\rm rank}(T_X) \over \varphi(m_X)}$.

\end{notation}

Corollary  \ref{characterization} characterizes the possible pairs $(n_X,r_X)$. 
This suggests the following question : 

\begin{question} What are the possiblities for $(n_X,m_X,r_X,t_X)$ ?

\end{question}

Huybrechts asked for explicit examples of $K3$ surfaces with $m_X > n_X$ (see \cite{H}, Remark 15.1.13). The following 
construction is due to Brandhorst and Elkies \cite{BE} :

\begin{example}\label{BE example} Let $X$ be the $K3$ surface constructed by Brandhorst and Elkies in \cite{BE}. 
This surface has an automorphism $a : X \to X$ such that the characteristic polynomial of $a^*|T_X$ is $\Phi_{14}$ (see \cite{BE},
\S 3), hence $m_X = 14$. Since $14 \not \in A \cup B$, we have $n_X \not = 14$, therefore $n_X < m_X$. 

\medskip This implies that $n_X = 7, 2$ or $1$. The results of Vorontsov \cite{V} and Oguiso-Zhang \cite{OZ}  imply that up to isomorphism, there 
exists a unique projective $K3$ surface $Y$ such that ${\rm rank}(T_Y) = 6$ and 
$n_Y = 7$ (note that this also follows from Proposition \ref{VK}); the discriminant of $S_Y$ is equal to $7$ (see \cite{OZ}, Lemma 1.3). 
For the $K3$ surface $X$ constructed by
Brandhorst and Elkies, the discriminant of $S_X$ is $7.13^2$, hence $X$ and $Y$ are not isomorphic; therefore $n_X \not = 7$.
One can show that $n_X = 1$.

\end{example}



\bigskip
\bigskip
Eva Bayer--Fluckiger 

EPFL-FSB-MATH

Station 8

1015 Lausanne, Switzerland

\medskip

eva.bayer@epfl.ch

\end{document}